\documentclass[11pt]{article}
\usepackage[active]{srcltx}
\usepackage{amssymb,amsmath}
\usepackage{color}
\usepackage[unicode]{hyperref}
\usepackage{txfonts}
\usepackage{graphicx}

\hypersetup{
    colorlinks = true,%
    citecolor = [rgb]{0.0,0.0,0.9},
    filecolor=black,%
    linkcolor = [rgb]{0.65,0.0,0.0},%
    anchorcolor = red,
    pagecolor = red,
    urlcolor= [rgb]{0.65,0.0,0.0}
}

\newtheorem{example}{ Example}[section]
\newtheorem{proposition}{Proposition}[section]
\newtheorem{theorem}{Theorem}[section]
\newtheorem{lemma}{Lemma}[section]

\newtheorem{definition}{Definition}[section]
\newtheorem{remark}{Remark}[section]

\numberwithin{equation}{section}

\def\l{\lambda}
\def\N{\mathbb{N}}
\def\H{\mathbb{H}}
\def\Hn{{\mathbb{H}^n}}

\def\R{\mathbb{R}}
\def\l{\lambda}

\def\t{\tau}
\def\a{\alpha}
\def\psn{\par\medskip\noindent}
\def\S{\mathbb{S}}

\begin{document}

\title{The engulfing property for sections of convex functions in the Heisenberg group and the associated quasi--metric}
%\title{Quasi--metrics  in the Heisenberg group via the engulfing property for convex functions}
\author{
 A. Calogero\thanks{
Dipartimento di Matematica e Applicazioni, Universit\`a degli Studi
di Milano-Bicocca, Via R. Cozzi 55, 20125 Milano, Italy ({\tt
andrea.calogero@unimib.it}, $^+$corresponding author: {\tt
rita.pini@unimib.it})},
\quad R. Pini$^*$$^+$
 }

\date{}
\maketitle

\begin{abstract}
\noindent  In this paper we investigate the property of engulfing for $H$-convex functions defined on the Heisenberg group $\Hn$. Starting from the horizontal sections introduced by Capogna and Maldonado in \cite{CapMal2006}, we consider a new notion of section, called $\Hn$-section, as well as a new condition of engulfing associated to the $\Hn$-sections, for an $H$-convex function defined in $\H^n.$ These sections, that arise as suitable unions of horizontal sections, are dimensionally larger; as a matter of fact, the $\Hn$-sections, with their engulfing property, will lead to the definition of a pseudo-metric in $\Hn$ in a way similar to Aimar, Forzani and Toledano in the Euclidean case (\cite{AiFoTo1998}). A key role is played by the property of round $H$-sections for an $H$-convex function, and by its connection with the engulfing properties.
%
%
%
%
%
%
%Similarly to the Euclidean setting, we introduce the notion of round $H$-sections for an $H$-convex function, and study its %connection  with the engulfing properties.
\end{abstract}



\section{Introduction}

Given a convex function $u:\R^n\to\R$, for every $x_0\in \R^n,$ $p\in \partial u(x_0),$ and $s>0,$  we will denote by $S_{u}(x_0,p,s)$ the section of $u$ at $x_0$ with height $s$, defined as follows
 \begin{equation}\label{def section R}
S_u(x_0,p,s)=\left\{x\in \R^n:\ u(x)- u(x_0)-p\cdot(x-x_0)<s\right\};
\end{equation}
in case $u$ is differentiable at $x_0,$ we will denote  the section  by $S_u(x_0,s),$ for short. The related notion of engulfing for convex functions, or, equivalently, for their sections, is essentially a geometric property, and it is based on a regular mutual behaviour of the sections of the function.
We say that a convex function $u$ satisfies
the \emph{engulfing property}
(shortly, $u \in E(\R^n,K)$) if there exists $K>1$ such that for any $x\in\R^n,$ $p\in \partial u(x),$ and
$s > 0$, if $y\in S_u(x,p,s)$, then $S_u(x,p,s)\subset S_u(y,q,Ks),$ for every $q\in \partial u(y)$.
%\end{definition}

The functions $u$ in the class $E(\R^n,K)$ have been studied in connection with the solution to the  Monge-Amp\`{e}re equation $\mathrm{det} \,D^2u=\mu,$ where $\mu$ is a Borel measure on $\R^n.$  %Several conditions on such functions $u$ have been proposed in order to preserve the harmony between measure theory (and in particular the Monge-Amp\`{e}re measure $\mu_u$), and the shape of the sections with their induced geometry.
In this framework, a $\mathcal{C}^{1,\beta}$-estimate for the strictly convex, generalized solutions to the Monge-Amp\`{e}re equation was proved by  Caffarelli (\cite{Caf1991, Caf1992}), under the assumption that the measure $\mu$ satisfies a suitable doubling property (see the exhaustive book by Guti\'errez \cite{Gu2001}).
This doubling property is actually equivalent to the geometric property of engulfing for the solution $u.$

Another issue is related to the properties enjoyed by the family of sections $\{S_u(x,s)\}_{\{x\in \R^n,\, s>0\}},$ in case $u$ is a convex differentiable  function in  $E(\R^n, K).$ In \cite{AiFoTo1998}, it is shown that, in this case, one can define a quasi-metric $d$ on $\R^n$ as follows:
\begin{equation}\label{intro dista}
d(x,y):=\inf \left\{s>0:\, x\in S_{u}(y,s),\,y\in S_{u}(x,s)\right\}.
\end{equation}
In addition, if $B_d(x,r)$ is a $d$-ball of center $x$ and radius $r,$ then
\begin{equation}\label{intro inclusione sezioni}
S_u \left(x,\frac{r}{2K}\right)\subset B_d(x,r)\subset S_u(x,r).
\end{equation}
In the archetypal case $u(x)=\|x\|,$ with $x\in\R^n$, one has $S_u(x,s)=B^{\R^n}(x,\sqrt{s}),$ and hence the family of sections of $u$ consists of the usual balls in $\R^n$.

In the case of convex functions defined in a Carnot group $\textbf{G}$, in \cite{CapMal2006} Capogna and Maldonado introduced some appropriate geometric objects, that can be considered as the sub-Riemmannian analogue of the classical sections, as well as a naturally related notion of horizontal engulfing.
Given a horizontally convex function $\varphi: \textbf{G}\to \R,$ $\xi_0\in \textbf{G},$ $p\in \R^{m_1},$ $s>0,$ the section $S_u^H(\xi_0,p,s)$ ($H$-sections, from now on, where $H$ stands for \emph{horizontal}) is defined as follows:
\begin{equation}\label{intro H section}
S_\varphi^H(\xi_0,p,s):=\{\xi_0\circ \exp v:\, v\in V_1, \,\varphi(\xi_0\circ \exp v)-\varphi(\xi_0)- v\cdot p <s\},
\end{equation}
where $V_1\cong \R^{m_1}$ is the first layer of the stratification of the Lie algebra of $\textbf{G}$;
in case $\varphi$ is horizontally differentiable at $\xi_0,$ we will denote such $H$-section  by $S_\varphi^H(\xi_0,s),$ for short.
The mentioned authors
say that a horizontal convex and differentiable function $\varphi$ satisfies the \emph{engulfing property} if there exists $K>1$ such that, for every $\xi,\, \xi'\in\textbf{G}$  and $s>0$,  if $\xi'\in S_\varphi^H(\xi,s)$,  then  $\xi\in  S_\varphi^H(\xi',Ks)$. Let us stress that the definition of $H$-section in (\ref{intro H section}) and the notion of engulfing  are affected by the sub-Riemannian structure exactly as the notion of horizontal convexity; more precisely, they rely upon the behaviour of the function on the horizontal lines and planes.
In \cite{CapMal2006} it is proved that the horizontal derivatives of a strictly convex and everywhere differentiable function on a Carnot group, satisfying this horizontal version of engulfing, belong to the Folland-Stein class $\Gamma^{1+1/K},$ i.e., the horizontal derivatives $X_i\varphi$ are $1/K$-H\"{o}lder continuous with respect to any left-invariant and homogeneous pseudo-norm in the group. The key point in their argument is a reduction of the general discussion to the one-dimensional case.  As a matter of fact, the topological dimension of the $H$-sections in (\ref{intro H section}) is the dimension of the first layer of the stratification of the Lie algebra of the group, and this prevents from building a pseudo-metric as in (\ref{intro dista}) starting from the family of sections associated to every point of the group.

In this paper we focus on horizontal convex functions $\varphi$ ($H$-convex functions) on the Heisenberg group $\H^n,$ that is the simplest Carnot group of step 2. Our main purpose is to overcome the dimensional gap between the $H$-sections defined in \cite{CapMal2006}, and the balls related to any pseudo-distance in $\H^n,$ by introducing and studying a different notion of section.
Our idea takes inspiration from the notion of $H$-section in (\ref{intro H section}), together with the property that any pair of points in $\H^n$ can be joined by at most three consecutive horizontal segments. These facts lead us to define full-dimensional sections that arise as a sort of composition in three steps of \lq\lq thin\rq\rq\,  $H$-sections.  These new objects will be called $\Hn$-section, and  will be denoted by
$\S_\varphi^\Hn(\xi_0,p,s)$ (for the precise definition of $\S_\varphi^\Hn(\xi_0,p,s),$ see Definition \ref{def section Heis}). For these $\Hn$-sections, we introduce the following engulfing condition:
\begin{definition}\label{intro defin engulfing cicciotta}
Let  $\varphi :\Hn \to \R$ be an $H$-convex function.
We say that
$\varphi$
satisfies
the \textbf{engulfing property} $E(\Hn,K)$ if there exists $K > 1$
such that for any $\xi\in\Hn,\ p\in \partial_H \varphi(\xi)$ and
$s > 0$, if $\xi'\in \S_\varphi^\Hn(\xi,p,s)$, then $$\S_\varphi^\Hn(\xi,p,s)\subset  \S_\varphi^\Hn(\xi',q,Ks),$$ for every $q\in \partial_H \varphi(\xi')$.
\end{definition}
It is obvious that a function which satisfies this engulfing property $E(\Hn,K)$, satisfies the engulfing property introduced by Capogna and Maldonado.

The study of this new notion of engulfing for $\Hn$-sections of full-dimension requires a mix of tools and properties inherited by the Euclidean case $\R^n$, both for the simplest case $n=1,$ and for the knotty case $n>1$. Following the idea in \cite{KoMa2005} and, in particular, the equivalence between iii. and iv in Theorem \ref{teo kovalev maldonado} below, we introduce and study a horizontal notion of round sections for the $H$-sections (see Definition \ref{def H round section}). We prove that every $H$-convex function with round $H$-sections satisfies the engulfing property $E(\Hn,K)$ in Definition \ref{intro defin engulfing cicciotta}.

Let us summarize our results as follows:
\begin{theorem}\label{intro Teo}
 Let $\varphi :\Hn \to \R$ be an $H$-convex function with round $H$-sections, then
 \begin{itemize}
 \item[i.] $\varphi$ satisfies the engulfing property $E(\Hn,K)$; consequently, in the class of $H$-convex functions with round $H$-sections, the engulfing for $H$-sections and the engulfing for $\Hn$-sections are equivalent properties;

 \item[ii.] the function $d_\varphi:\Hn\times\Hn\to [0,+\infty)$ defined by
$$
d_\varphi(\xi,\xi')=\inf\left\{s>0:\ \xi\in \S_\varphi^\Hn(\xi',s),\  \xi'\in \S_\varphi^\Hn(\xi,s)\right\}
$$
is a quasi-metric in $\Hn$; moreover, for the $d_\varphi$-balls, an $\Hn$-version of the inclusions in (\ref{intro inclusione sezioni}) holds true (see (\ref{dim finale balle}) below).
\end{itemize}
 \end{theorem}
Here, the archetypal example in $\H$ of the $H$-convex function $\varphi(x,y,t)=x^2+y^2$ gives  $S_\varphi^\H(\xi,s)=\widetilde B(\xi,\sqrt{s}),$ that is, the family of $\H$-sections of $\varphi$ consists of the $\widetilde B$-balls of a left-invariant and homogeneous distance $\widetilde d$ (see (\ref{balls tilde}) and
Example \ref{esempio}).

The property of round $H$-sections is actually stronger that the horizontal engulfing; we are able to provide an example of an $H$-convex function which satisfies the horizontal engulfing property  but has not round $H$-sections, and this phenomenon appears also in the Euclidean case, if $n>1$. Nevertheless, the main issue of the result above relies upon the dimensional gap between the assumptions, where a purely horizontal property is required, and the final result, where full-dimensional sets are involved.

The paper is organized as follows. In Section 2 we recall some results related to the engulfing property for a function defined in $\R^n,$ together with the structure of $\H^n$ and the notion of horizontal convexity. In Section 3 we introduce the $H$-sections, and we show that  round $H$-sections and controlled $H$-slope are equivalent property for these $H$-sections (see Theorem \ref{theorem mM se e solo se round}). In Section 4 we characterize the functions with the engulfing property $E(H,K),$ and prove that the two properties introduced in the previous section are sufficient conditions for a function to be in $E(H,K)$. In Section 5 we move to the notion of $\Hn$-sections and the related engulfing property as in Definition \ref{intro defin engulfing cicciotta}, and we prove Theorem \ref{intro Teo} i. In Section 6 we prove  Theorem \ref{intro Teo} ii. and provide a concrete example. In the final section we list some open questions.

\section{Preliminary notions and results}
In the paper, we will deal with $H$-convex functions defined on the Heisenberg group $\H^n.$ As we will see later, the notion of $H$-convexity requires that, for every point $\xi\in \H^n,$ one looks at the behaviour of the function under two points of view. The first one is one-dimensional, since the restriction of the function to any horizontal line $\{\xi\circ \exp tv\}_{t\in \R},$ with $v\in V_1,$ is an ordinary convex function; the second one is $2n$-dimensional, according to the fact that $v\in V_1\cong \R^{2n},$ or, equivalently, the horizontal lines through $\xi$ span the $2n$-dimensional horizontal plane $H_\xi.$
For these reasons, the first part of this section will be devoted to some results related to the engulfing property of convex functions $u:\R^n\to\R$, both in the case $n=1,$ and in the case $n\ge 2.$ In the second part we will recall the notion of $H$-convexity, together with some related results, for functions defined in the Heisenberg group $\Hn$.

\subsection{The engulfing property for convex functions in $\R^n$}

Let us concentrate, first, on the one-dimensional case, i.e. $n=1$.
The following characterization holds (see Theorem 2 in \cite{FoMa2004}, Theorem 5.1 in \cite{CUFoMa2008}):
\begin{theorem}\label{th:softR}
Let $u:\R\to\R$ be a strictly convex and differentiable function. The following are equivalent:
\begin{itemize}
\item[i.] $u\in E(\R,K)$, for some $K>1$;
\item[ii.] there exists a constant $K'>1$ such that, if $x,y\in \R$ and $s>0$ verify $x\in S_u(y,s),$ then $y\in S_u(x,K's);$
    \item[iii.] there exists a constant $K''> 1$ such that, for any $x,\ y\in\R,$
\begin{equation}\label{carat 3 in R}
\begin{array}{l}
\displaystyle \frac{K''+1}{K''}\left(u(y)-u(x)-u'(x)(y-x)\right)\le (u'(x)-u' (y))(x-y)\\
\qquad\qquad\qquad\qquad\qquad\le (K''+1)\left(u(y)-u(x)-u'(x)(y-x)\right).
\end{array}
\end{equation}
\end{itemize}
\end{theorem}
As a matter of fact, the assumption of differentiability in the theorem above can be removed, as proved in  \cite{CaPi2020}:
\begin{theorem}\label{teo differenziabilita}
Let $u:\R\to\R$ be a convex function, with bounded sections, satisfying the engulfing property. Then, $u$ is strictly convex and is in ${\mathcal C}^1(\R)$.
\end{theorem}

Given a strictly convex differentiable function $u:\R\to \R,$ one can consider the associated Monge-Amp\`{e}re measure $\mu_{u}$ defined on any Borel set $A\subset \R$ by
$$
\mu_{u}(A)=|u'(A)|,
$$
where $|\cdot |$ denotes the Lebesgue measure. We say that
the measure $\mu_{u}$ has the \emph{(DC)-doubling property} if there exist constants $\alpha\in (0,1)$ and $C>1$ such that
\begin{equation}\label{DC}
\mu_{u}(S_{u}(x,s))\le C\mu_{u}(\alpha S_{u}(x,s)),
\end{equation}
for every section $S_{u}(x,s)$ (here $\alpha S_{u}(x,s)$ is the open convex set obtained by $\alpha$-contraction of $S_{u}(x,s)$ with respect to its center of mass).
In \cite{GuHu2000} and \cite{FoMa2002} it was shown that the (DC)-doubling property of the measure $\mu_{u}$ is equivalent to the engulfing property for the function $u;$ in particular, given $u$ in $E(\R,K),$  the constants $\alpha$ and $C$ in (\ref{DC}) depend only on $K.$
A Radom measure $\mu$ is \emph{doubling} if and only if there exists a constant $A$ such that
\begin{equation}\label{doubling}
\frac{1}{A}\le\frac{\mu(Q_1)}{\mu(Q_2)}\le A,
\end{equation}
for any congruent cubes $Q_1$ and $Q_2$ with nonempty intersection (see, for example, \cite{KoMaWu2007}). We recall that two subsets of $\R$ are called \emph{congruent} in there exists an isometry of $\R$ that maps one of them onto the other.
Since every open and bounded interval in $\R$ is a particular section for $u$, the (DC)-doubling property of $\mu_{u}$ is trivially equivalent to the fact that $\mu_{u}$ is a doubling measure. In particular, the constant $A$ depends only on $K.$
Now, noticing that $\mu_{u}((x,x+r))=u'(x+r)-u'(x)$, by  (\ref{carat 3 in R}) we obtain
$$
\frac{K''+1}{K''}\left(u(x+r)-u(x)-u'(x)r\right)\le r\mu_{u}((x,x+r)))\le (K''+1)\left(u(x+r)-u(x)-u'(x)r\right).
$$
These arguments show the central role of the function $(x,r)\mapsto u(x+r)-u(x)-u'(x)r$ in our paper.
More precisely in \cite{CUFoMa2008} (see Theorem 5.5) the authors prove the following:
\begin{theorem}
Let $u:\R\to\R$ be a strictly convex and differentiable function. Then $u\in E(\R,K)$ if and only if
there exist two constants $A_1>1$ and $A_2>1,$ both of them depending on $K_,$ such that
\begin{eqnarray}
&&\frac{1}{A_1}\le \frac{u(x+r)-u(x)-u'(x)r}{u(x-r)-u(x)+u'(x)r}\le A_1,\qquad\forall x\in\R,\ r>0;\label{symmetric}\\
&&\frac{1}{A_2}\le \frac{u(x+2r)-u(x)-u'(x)2r}{u(x+r)-u(x)-u'(x)r}\le A_2,\qquad\forall x\in\R,\ r>0.\label{delta2}
\end{eqnarray}
\end{theorem}
Condition (\ref{symmetric}) says that $u$ is essentially symmetric around every point, and condition (\ref{delta2}) says that it satisfies the so-called $\Delta_2$ condition at each point in $\R$.

Hence, the behaviour of the measure $\mu_{u}$ is related to the
functions $m_u,\ M_u:\R\times\R^+\to\R^+$ defined by
\begin{equation}\label{eq:M_and_m}
\begin{aligned}
m_u(x,r)&:=\min_{\{z:\ |z-x|=r\}} \left(u(z)-u(x)-u'(x)(z-x)\right)\\
M_u(x,r)&:=\max_{\{z:\ |z-x|=r\}} \left(u(z)-u(x)-u'(x)(z-x)\right),
\end{aligned}
\end{equation}
for every $x\in\R,\ r\in\R^+$. These functions will be naturally extended to the $n$-dimensional case and in $\Hn$, and will play a crucial role in the investigation of the engulfing for $H$-convex functions.

For every fixed $x\in\R,$ denote by $u_x$ the function
\begin{equation}\label{ux}
s\mapsto u_x(s)=u(x+s)-u(x)-u'(x)s.
\end{equation}
Then, $M_u(x,r)\in\{u_x(\pm r)\},$ and $M_u(x,2r)\in\{u_x(\pm 2r)\}.$ Let
us suppose, for instance, that the following equalities hold true:
 $$ M_u(x,2r)=u_x(2r),\quad m_u(x,2r)=u_x(-2r),\quad
 M_u(x,r)=u_x(r), \quad  m_u(x,r)=u_x(-r).
 $$
 Then,
 by (\ref{symmetric}) and (\ref{delta2}), we obtain
\begin{eqnarray*}
&&M_u(x,2r)=u_x(2r)\le A_2 u_x(r)=A_2 M_u(x,r),\\
&&m_u(x,2r)=u_x(-2r)\le A_1 u_x(2r)\le A_1A_2 u_x(r)\le A_1^2A_2 u_x(-r)=A_1^2A_2m_u(x,r),\\
&&M_u(x,r)=u_x(r)\le A_1 u_x(-r)=A_1 m_u(x,r).
 \end{eqnarray*}
The other possible combinations can be treated similarly, and we obtain
the following  fundamental estimates:
 \begin{remark}\label{remark}{\rm
 Let  $u\in E(\R,K)$ be a strictly convex and differentiable function. Then,
  \begin{eqnarray}
&M_u(x,2r)\le B_1 M_u(x,r),\qquad &\forall x\in\R,\ r\ge 0\label{B0}\\
&m_u(x,2r)\le B_2 m_u(x,r),\qquad &\forall x\in\R,\ r\ge 0\label{B4}\\
&M_u(x,r)\le B_3 m_u(x,r),\qquad &\forall x\in\R,\ r\ge 0\label{B1}
\end{eqnarray}
where $B_1,\ B_2$ and $B_3$ depend only on $K$ (and $B_i>1$).}
 \end{remark}
It is worthwhile to note that inequality (\ref{B1}) is false if $n\ge 2,$ despite the engulfing property holds; the function in (\ref{ex Wang}), due to Wang, will provide a counterexample to this phenomenon.

The next result provides another estimate for the function $m_u$:
\begin{proposition}\label{prop new}
Let $u\in E(\R,K)$ be a convex function with bounded sections. Then,
 \begin{equation}\label{B2}
B_4 m_u(x,r)\le m_u(x,2r),\qquad \forall x\in\R,\ r\ge 0,
\end{equation}
with $B_4>1$ which depends only on $K$.
\end{proposition}
\psn
\noindent \textbf{Proof.}
Let us fix $x\in\R$. The function $u_x$ defined in (\ref{ux}) is strictly convex and differentiable (see \cite{CaPi2020}), and belongs to $E(\R,K)$; moreover,
$$
\frac{K''+1}{K''}u_x(y)\le u_x'(y)y,\qquad\forall y\in\R,
$$
where $K''$ depends only on $K$  (for all the details, see Theorem 4 and its proof in \cite{FoMa2004}). Hence, for every fixed $r>0$, the Gronwall inequality gives
$$
u_x(|y|)\ge u_x(\alpha)\left(\frac{|y|}{\alpha}\right)^{\frac{K''+1}{K''}},\qquad\forall |y|\ge r.
$$
Therefore, we obtain that $u_x(\pm 2r)\ge 2^{\frac{K''+1}{K''}}u_x(\pm r).$
Let $m_u(x,r)=u_x(r).$ Then, $B_4 m_u(x,r)\le u_x(2r).$ Suppose that $B_4m_u(x,r)>u_x(-r).$ In this case, $u_x(-2r)\ge B_4 u_x(-r),$ and thus $B_4u_x(-r)<B_4 m_u(x,r),$ a contradiction. Then, (\ref{B2}) follows.
\hfill$\square$\psn

Let us now move to the case $n\ge 2.$ Given a differentiable function $u:\R^n\to \R,$ as in the one-dimensional case \eqref{eq:M_and_m}, the functions
$m_u,\ M_u:\R^n\times\R^+\to\R^+$ are defined by
\begin{eqnarray*}
m_u(x,r)&=&\min_{\{z:\ \|z-x\|=r\}} \left(u(z)-u(x)-\nabla u(x)\cdot(z-x)\right)\\
M_u(x,r)&=&\max_{\{z:\ \|z-x\|=r\}} \left(u(z)-u(x)-\nabla u(x)\cdot(z-x)\right),
\end{eqnarray*}
for every $x\in\R^n,\ r\in\R^+$.

Let us recall the following property, that will be critical when dealing with the engulfing in $\H^n.$
\begin{definition}\label{def round sections}(see Definition 2.1 in \cite{KoMa2005}) Let $u:\R^n\to \R$ be a convex function. We say that $u$ has \emph{round sections} if there exists a constant $\tau\in (0,1)$ with the following property: for every $x\in \R^n,$ $p\in \partial u(x),$ and $s>0,$ there is $R>0$ such that
$$
B(x,\tau R)\subset S_u(x,p,s)\subset B(x,R).
$$
\end{definition}
In \cite{KoMa2005} (see Theorem \ref{teo kovalev maldonado} below) it is proved that a convex function $u:\R^n\to \R$ has round sections if and only if $u$ is differentiable, but not affine, and has \emph{controlled slope}, i.e., there exists a constant $H\ge 1$ such that
\begin{equation}\label{def controlled slope}
M_u(x,r)\le H m_u(x,r),\qquad \forall x\in\R^n,\ r\ge 0.
\end{equation}
This equivalence is quantitative, in the sense that the constants involved in each statement depend only on each other and $n,$ but not on $u.$
Furthermore, if $u:\R^n\to \R$ satisfies one of the two equivalent conditions above, then $u\in E(\R^n,K),$ for a suitable $K>1$ (see Theorem 3.9 in \cite{KoMa2005}).
Let us finally notice that condition (\ref{def controlled slope}) is the $n$-dimensional version of condition (\ref{B1}): in the case $n\ge 2$, hence, the controlled slope for a function, or, equivalently, the property of round sections, is only a sufficient condition for a function to have the engulfing property.

\subsection{Convexity in the Heisenberg group $\Hn$}

The Heisenberg group $\Hn$ is the simplest Carnot group of step 2. We will recall some of the notions and
background results used in the sequel. We will focus only on those geometric
aspects that are relevant to our paper.  For a general overview on the subject, we
refer to \cite{BoLaUg2007} and \cite{CaDaPaTy2007}.

The  Lie algebra  $\frak{h}$ of $\mathbb H^n$ admits a
stratification  $ \frak{h}=V_1\oplus V_2$ with $V_1=\texttt{\rm
span}\{X_{i},\, Y_i;\ 1\le i\le n\}$ being the first layer of the so-called horizontal vector fields, and $V_2=\texttt{\rm span}\{T\}$
being the second layer which is one-dimensional. We assume
$[X_i,Y_i]=-4T$ and the remaining commutators of basis vectors
vanish. The exponential map $\exp:\frak{h}\to \mathbb H^n$ is
defined in  the usual way. By these commutator rules we obtain,
using the Baker-Campbell-Hausdorff formula, that $\Hn$ can be identified with
$\R^{n}\times \R^{n}\times \R$  endowed
with the non-commutative group law given by
$$
\xi\circ \xi'=(x,y,t)\circ (x',y',t')=(x+x',y+y',t+t'+2 (x'\cdot y- x\cdot y')),
$$
where $x,y,x'$ and $y'$ are in $\R^n$, $t$ and $t'$ in $\R$, and where $'\cdot '$ is the inner product in $\R^n$. Let us denote by $e$ the neutral element in $\Hn.$
Transporting the basis vectors of $V_1$ from the origin to an
arbitrary point of the group by a left-translation, we obtain a
system of left-invariant vector fields written as first order
differential operators as follows
\begin{equation} \label{vector fields}
     X_j=\partial_{x_j}+2y_j \partial_t,\qquad
     Y_j=\partial_{y_j}-2x_j\partial_t,\qquad j=1,...,n.
\end{equation}

Via the exponential map $\exp:\mathfrak{h}\to\H$ we identify the
vector $\sum_{i=1}^n(\alpha_i X_i+\beta_i Y_i)+\gamma T$ in
$\mathfrak{h}$ with the point $(\alpha_1,\ldots,\alpha_n,
\beta_1,\ldots,\beta_n, \gamma)$ in $\Hn;$ the inverse $\xi:
\Hn\to\mathfrak{h}$ of the exponential map has the unique
decomposition $\xi=(\xi_1,\xi_2),$ with $\xi_i:\Hn\to V_i,$ and  we
identify $V_1$ with $\R^{2n}$ when needed.

For every positive $\lambda,$ the non-isotropic Heisenberg
dilation $\delta_{\lambda}: \H^{n}\to \H^{n}$ is
defined by $\delta_{\lambda}(x,y,t) = (\lambda x,\lambda y,
\lambda^{2}t)$. Let $N(x,y,t)=((\|x\|^2+\|y\|^2)^2+t^2)^\frac{1}{4}$ be the gauge
norm in $\mathbb H^n$. The function $d_g:\H^n\times \H^n\to [0,+\infty)$ defined by
$$
d_{g}(\xi,\xi'):=N((\xi')^{-1}\circ \xi)
$$
satisfies the triangle inequality, thereby defining a metric on $\Hn$: this
metric is the so-called Kor\'anyi-Cygan metric which is
left-invariant and homogeneous, i.e. $d_g(\delta_\lambda(\xi),\delta_\lambda(\xi'))=\lambda d_g(\xi,\xi')$ for every $\lambda>0,$ $\xi,\xi'\in\Hn$. We will set $d_g(e,\xi)=\|\xi\|_g$ for every $\xi\in\H^n$. The Kor\'anyi-Cygan  ball of center $\xi_0\in
\mathbb H^n$ and radius $r>0$ is given by
$B_{g}(\xi_0,r)=\{\xi\in \Hn:\ d_g(\xi_0,\xi)\le r\}.$

The horizontal structure relies on the notion of horizontal
plane. Given $\xi_0\in\Hn$, the {\emph{horizontal plane}}
$H_{\xi_0}$ associated to $\xi_0=(x_0,y_0,t_0)$ is the plane in
$\Hn$ defined by
$$
H_{\xi_0}:=\left\{\xi=(x,y,t)\in\Hn:\ t=t_0+2( y_0\cdot x- x_0\cdot y)\right\}.
$$
This is the plane spanned by the horizontal vector fields $\{X_i,\
Y_i\}_i$ at the point $\xi_0$; note that  $\xi'\in H_{\xi}$ if
and only if $\xi\in H_{\xi'}$. A \emph{horizontal segment} is a convex subset of a
\emph{horizontal line}, which is a line lying on a horizontal plane $H_\xi$ and passing though  the point $\xi\in\Hn$;
if $\xi'\in H_{\xi}$, with $\xi'\not=\xi$, then $H_\xi\cap H_{\xi'}$ is a horizontal line.

Let $\Omega\subset \mathbb H^n$ be an open set.   The main idea of
the analysis in the Heisenberg group is that the regularity
properties of functions defined in $\H^n$ can be
expressed in terms only of the horizontal vector fields \eqref{vector
fields}. In particular, the appropriate notion of gradient for a
function is the so-called {\it horizontal gradient}, which is
defined as the $2n$-vector $ \nabla_H\varphi(\xi)=
\left(X_1\varphi(\xi),...,X_n\varphi(\xi), Y_1\varphi(\xi),...,Y_n\varphi(\xi)\right)$ for a
function $\varphi\in\Gamma^1(\Omega)$. % at $\xi\in\Omega$ is the
Here, $\Gamma^k(\Omega)$ denotes the Folland--Stein space of
functions having continuous derivatives up to order $k$ with respect
to the vector fields $X_i$ and $Y_i,$ $i\in \{1,...,n\}$.
We say that $\varphi : \Omega\to\R$  is $H$-\emph{differentiable} at $\xi$, if there exists
a mapping $D_H \varphi :\Hn\to\R$  which is $H$-linear, i.e. $D_H \varphi (x,y,t)=D_H \varphi (x,y,0)$ for every $(x,y,t)\in\Hn$, such that $\varphi(\xi\circ\xi') = \varphi(\xi)+ D_H \varphi(\xi')+o(\|\xi'\|_g)$; the  vector associated to $D_H \varphi$ with respect to
the fixed scalar product is the horizontal gradient $ \nabla_H \varphi(\xi)$.

For general non-smooth functions $\varphi:\Omega\to \mathbb R,$ the
{\it horizontal subdifferential} $\partial _H \varphi(\xi_0)$ of $\varphi$ at
$\xi_0\in \Omega$ is
 given by
$$
\partial _H \varphi(\xi_0)=\left\{p\in \mathbb R^{2n}:\varphi(\xi)\geq \varphi(\xi_0)+p\cdot
({\texttt{\rm Pr}_1}(\xi)-{\texttt{\rm Pr}_1}(\xi_0)), \ \forall
\xi\in \Omega\cap H_{\xi_0}\right\},
$$
where ${\texttt{\rm
Pr}_1}:\H^n\to\R^{2n}$ is the projection defined by ${\texttt{\rm
Pr}_1}(\xi)={\texttt{\rm Pr}_1}(x,y,t)=(x,y)$.
It is easy to see that if $\varphi\in \Gamma^1(\Omega)$ and
$\partial_H\varphi(\xi)\neq \emptyset$, then $\partial_H\varphi(\xi)=\{\nabla_H
\varphi(\xi)\}$.
  A function $\varphi:\Omega\to \mathbb R$ is called $H-${\it
subdifferentiable on} $\Omega$ if $\partial_H \varphi(\xi)\neq \emptyset$
for every $\xi\in \Omega.$

A central object of study within this paper is provided by the $H$-convex functions. First of all, we recall that
a set $\Omega\subset \mathbb H^n$ is said to be {\it horizontally convex} ($H$-{\it convex}) if, for every
$\xi_1,\xi_2\in \Omega,$ with $\xi_1\in H_{\xi_2}$ and
$\lambda\in [0,1]$, we have $\xi_1\circ \delta_\lambda
(\xi_1^{-1}\circ \xi_2)\in \Omega$. It is clear that if
$\Omega$ is convex (i.e. it is convex in the $\R^{2n+1}$-sense),
then it is also $H$-convex.
Given a function $\varphi:\Omega\to \mathbb R,$ where $\Omega$ is $H$-convex, there are several equivalent ways to define the concept of $H$-convexity for $\varphi.$
The most intuitive one is to require  the classical convexity of the function when restricted to any horizontal line within $\Omega.$ The same definition can be rephrased by considering the group operation: the function
$\varphi:\Omega\to \mathbb R$ is said to be $H$-{\it convex}
if, for every $\xi_1,\xi_2\in  \Omega$ with $\xi_1\in
H_{\xi_2}$ and $\lambda\in [0,1]$, we have that
\begin{equation}\label{convex-def}
    \varphi(\xi_1\circ \delta_\lambda (\xi_1^{-1}\circ \xi_2))\leq
(1-\lambda)\varphi(\xi_1)+\lambda \varphi(\xi_2).
\end{equation}
 If the  strict inequality holds in (\ref{convex-def}), for every
$\xi_1\neq \xi_2$ and $\lambda\in (0,1)$,   then $\varphi$ is said to be {\it
strictly} $H$-{\it convex}.
$H$-convex functions have been extensively studied in the last few years;
their characterizations, as well as their regularity
properties, like their continuity, for instance, will come into play through the paper, and we
  refer to
\cite{BaRi2003,CaPi2011,DaGaNh2003,MaSc2014}. Let us recall, in particular, that $\varphi:\H^n\to\R$ is $H$-convex if and only if $\varphi$ is $H$-subdifferentiable.

\section{$H$-convex functions with round $H$-sections and with controlled $H$-slope}

As already seen in the Introduction, a horizontal notion of section was given in \cite{CapMal2006} for functions defined on a general Carnot group $\textbf{G}.$ We will consider the particular case $\textbf{G}=\H^n.$

Let $\varphi:\Hn \to\R$ be an $H$-convex function, and let us fix $\xi_0\in\Hn,$ $p_0\in \partial_H \varphi(\xi_0),$  and $s>0.$ The \emph{$H$-section} of $\varphi$ at $\xi_0,$ $p_0,$ with height $s,$ is the set
\begin{equation}\label{def section H}
S_\varphi^H(\xi_0,p_0,s)=\{\xi\in H_{\xi_0}:\ \varphi(\xi)- \varphi(\xi_0)-p_0\cdot
({\rm Pr}_1(\xi)-{\rm Pr}_1(\xi_0))<s\}.
\end{equation}
If $\varphi$ is $H$-differentiable, then $\partial_H \varphi(\xi_0)=\{\nabla_H \varphi(\xi_0) \},$ and we simply write
$S_\varphi^H(\xi_0,s)$ for the corresponding $H$-section.
For every fixed $(\xi_0,p_0,s),$ the set $S_\varphi^H(\xi_0,p_0,s)$ is $H$-convex, and is contained in a horizontal plane; this dimensional gap between $H$-sections and open sets in $\Hn$ is a crucial difference with respect to the Euclidean case.

In this section we essentially introduce the notions of round $H$-sections (see Definition \ref{def H round section}) and controlled $H$-slope (see Definition \ref{def m M}), proving their equivalence (see Theorem \ref{theorem mM se e solo se round}). Let us emphasize that these two properties for an $H$-convex function are horizontal properties, i.e. they give information on the behaviour of the function only when restricted to the horizontal planes, exactly as the notion of $H$-section, $H$-convexity and $H$-subdifferential.

In the following of the paper, for every function $\varphi:\Hn \to\R,$ and for every $\xi_0\in \Hn$, $p_0\in\partial_H\varphi(\xi_0)$ and  $v_0\in V_1\setminus \{0\}$, we will consider the functions
$\varphi_{\xi_0,p_0}:\Hn\to\R$ and $\widehat\varphi_{\xi_0,v_0}:\R\to\R$ defined by
\begin{eqnarray}
&&\varphi_{\xi_0,p_0}(\xi)=\varphi(\xi)-\varphi(\xi_0)-p_0\cdot(\texttt{\rm Pr}_1(\xi)-\texttt{\rm Pr}_1(\xi_0)),\qquad \forall \xi\in \Hn,\label{def varphi xi}\\
&&
\widehat\varphi_{\xi_0,v_0}(\a)=\varphi(\xi_0\circ\exp (\a v_0)),\qquad \forall \a\in \R.\label{def widehat varphi xi0 v}
\end{eqnarray}
If $\varphi$ is $H$-differentiable, then we will set $\varphi_{\xi_0,\nabla_H\varphi(\xi_0)}=\varphi_{\xi_0}$.
The following result holds:
\begin{proposition}\label{proposizione SH limitate}
Let $\varphi :\Hn \to \R$ be a  strictly $H$-convex function. Then, all its $H$-sections are bounded sets.
\end{proposition}
\psn
\noindent \textbf{Proof.}  For every $\xi_0\in \Hn$ and $v\in V_1\setminus \{0\}$ let us consider
the function $\widehat\varphi_{\xi_0,v}$ as in (\ref{def widehat varphi xi0 v}).
By contradiction, let us suppose that there exists a sequence $\{(v_n,\alpha_n)\}_n$, with $v_n\in V_1,\ \|v_n\|=1,\ \alpha_n\to +\infty,$ such that $\xi_0\circ\exp(\alpha_n v_n)\in S^H_\varphi(\xi_0,p_0,s)$. Clearly, there exists a subsequence such that $v_n\to v_0\in V_1.$

Let us denote by $\a_0=\sup\left\{\a\ge 0:\ \xi_0\circ\exp(\a v_0)\in \overline{S^H_\varphi(\xi_0,p_0,s)}\right\}.$ If $\a_0=+\infty$, then the section $S_{\widehat\varphi_{\xi_0,v}}(0,s)$ of the function
$\widehat\varphi_{\xi_0,v}$ is unbounded; this is impossible, since $\widehat\varphi_{\xi_0,v}$ is strictly convex.
Let $s_0$ be finite, and let us consider the function
$\varphi_{\xi_0,p_0}$ in (\ref{def varphi xi});
the set $A=\{\xi\in\Hn:\ \varphi_{\xi_0,p_0}(\xi)\le s\}$ is  $H$-convex, since the function $\varphi_{\xi_0,p_0}$ is $H$-convex. Now, the previous arguments give
\begin{equation*}\left\{\xi\in H_{\xi_0}:\ \xi=\xi_0\circ\exp(\a v_n),\ 0\le \a\le \alpha_n\right\}\subset A,\; \forall n, \quad \mathrm{and}\quad
\xi'=\xi_0\circ\exp(\a_0 v_0)\in \partial A.
\end{equation*}
This contradicts Theorem 1.4 in \cite{ArCaMo2012}.
\hfill$\square$\psn

The next definition is related to a purely geometric property of the sections, and it will play a crucial role in the following of the paper.
\begin{definition}\label{def H round section}
We say that an $H$-convex function $\varphi : \Hn\to \R$ has \emph{round $H$-sections}  if there
exists a constant $K_0\in(0, 1)$  with the following property: for every $\xi\in \Hn$,
$p\in\partial_H \varphi(\xi)$ and $s>0,$ there exists $R>0$ such that
\begin{equation}\label{round section}
B_g(\xi, K_0R)\cap H_\xi\subset S^H_\varphi(\xi, p, s)\subset B_g(\xi,R)\cap H_\xi.
\end{equation}
\end{definition}
In particular, (\ref{round section}) implies that every $H$-section of a function with round $H$-sections is a bounded set. Clearly, Definition
\ref{def H round section} is the $\Hn$-version of Definition \ref{def round sections}; let us stress that it relies upon the subriemannian structure of $\Hn$ since, for every point $\xi,$ we restrict our attention only to the horizontal plane $H_\xi$.

\begin{remark}\label{rm:subdifferential} Let $\varphi:\Hn\to \R$ be $H$-convex, and consider the convex function $\widehat\varphi_{\xi_0,v}:\R\to \R$ defined by
(\ref{def widehat varphi xi0 v}). Then, if the nonempty convex set $\partial_H\varphi(\xi_0)$ is not a singleton, there exists $v\in V_1$ such that $\partial \widehat\varphi_{\xi_0, v}(0)$ is not a singleton. Indeed, suppose that $p+\l q\in \partial_H\varphi (\xi_0),$ for every $\l\in [0,1],$ with $q\neq 0.$ Then, by taking $v=q,$ we have that
$$
\widehat\varphi_{\xi_0,q}(\a)=\varphi(\xi_0\circ\exp (\a q))\ge \widehat\varphi_{\xi_0,q}(0)+\a (p\cdot q+\l\|q\|^2),\qquad \forall \l\in[0,1],\ \a\in\R.
$$
Hence $p\cdot q+\l\|q\|^2\in \partial \widehat\varphi_{\xi_0,q}(0)$ for every $\l\in [0,1].$
This implies that, if $\widehat\varphi_{\xi_0,v}$ is differentiable at $0$ for every $v\in V_1,$ then $\varphi$ is $H$-differentiable at $\xi_0.$
\end{remark}
In the previous remark and in the following result, the $H$-convexity plays a fundamental role in order to obtain some regularity properties of the function involved.

\begin{proposition}\label{round implica Hdiffenziabilita}
If $\varphi : \Hn\to \R$ is an $H$-convex function with round $H$-sections, then it is $H$-differentiable and strictly $H$-convex. Moreover, there exists a constant $C$ such that, for every $\xi_0\in \Hn$ and $v\in V_1,$ we have
\begin{equation}\label{pro 1}
\varphi_{\xi_0}(\xi_0\circ\exp(2v))\le C \varphi_{\xi_0}(\xi_0\circ\exp v),
\end{equation}
where the constant $C$ depends only on $K_0$ in  (\ref{round section}).
\end{proposition}

\noindent \textbf{Proof.} First of all note that, for every $\xi_0\in \Hn$ and $v\in V_1\setminus\{0\},$ the function $\widehat\varphi_{\xi_0,v}$ defined in (\ref{def widehat varphi xi0 v})  is convex, with round sections (with constant $K_0$). Therefore, Lemma 3.2 in \cite{KoMa2005} implies that it is differentiable and strictly convex. In particular, $\varphi$ is strictly $H$-convex.
Let us first show that $\varphi$ is $H$-differentiable at $\xi_0\in\Hn$. Since $\varphi$ is $H$-convex, this is equivalent to prove that the nonempty convex set $\partial_H\varphi(\xi_0)$ is a singleton (see Theorem 4.4, Prop. 5.1 in \cite{CaPi2011}, Theorem 1.4 in \cite{MaSc2014}). Suppose, by contradiction, that $\partial_H \varphi(\xi_0)$ is not a singleton; then, by Remark \ref{rm:subdifferential}, there exists $v\in V_1$ such that $\partial \widehat\varphi_{\xi_0,v}(0)$ is not a singleton. This
contradicts the fact that $\widehat\varphi_{\xi_0,v}(0)$ is differentiable.

Finally, taking into account that the function $\widehat\varphi_{\xi_0,v}$ is convex, differentiable and with round sections with constant $K_0$ , for every $\xi_0\in \Hn$ and $v\in V_1,$  again, by Lemma 3.2 in \cite{KoMa2005}, one has that there exists a constant $C$ depending only on $K_0$ such that
$$
\varphi_{\xi_0}(\xi_0\circ 2v)\le C \varphi_{\xi_0}(\xi_0\circ \exp v).
$$
\hfill$\square$\psn

In the sequel, given an $H$-differentiable function $\varphi : \Hn\to \R$, we will deal with the functions $m_\varphi^H,\ M_\varphi^H:\Hn\times\R^+\to\R^+$ that will take the place in $\H^n$ of the functions $m_u$ and $M_u$ in $\R^n.$ They are defined as follows:
\begin{eqnarray*}
m_\varphi^H(\xi,r)&:=&\min_{\{\xi'\in H_\xi:\ d_g(\xi,\xi')=r\}} \left(\varphi(\xi')-\varphi(\xi)-\nabla_H \varphi (\xi)\cdot(\texttt{\rm Pr}_1(\xi')-\texttt{\rm Pr}_1(\xi))\right)\\
M_\varphi^H(\xi,r)&:=&\max_{\{\xi'\in H_\xi:\ d_g(\xi,\xi')=r\}} \left(\varphi(\xi')-\varphi(\xi)-\nabla_H \varphi (\xi)\cdot(\texttt{\rm Pr}_1(\xi')-\texttt{\rm Pr}_1(\xi))\right),
 \end{eqnarray*}
 for every $\xi\in\Hn,\ r>0$.

A simple exercise shows that, if $\varphi :\Hn \to \R$ is an  $H$-differentiable and strictly $H$-convex function, then
for every $ \xi\in\Hn,$ and $r>0$,
\begin{equation}\label{inclusione con palla 2}
S^H_\varphi(\xi,   m_\varphi^H(\xi,r))\subset B_g(\xi,r)\cap H_{\xi}\subset S^H_\varphi(\xi,  M_\varphi^H(\xi,r)).
\end{equation}

The next definition, inherited from the corresponding one in $\R^n$ (see (\ref{def controlled slope})), pertains to the mutual behaviour of $m_\varphi^H$ and $M_\varphi^H,$ always from a horizontal point of view:
\begin{definition}\label{def m M}
We say that an $H$-convex function $\varphi : \Hn\to \R$ has \emph{controlled $H$-slope} (controlled horizontal slope) if $\varphi$ is $H$-differentiable, and there
exists a constant $K_1>0$ such that,
 for every $\xi\in \Hn$ and $r>0,$
\begin{equation}\label{H mM}
M_\varphi^H(\xi,r)\le K_1 m_\varphi^H (\xi,r).
\end{equation}
\end{definition}
Like in the Euclidean case (see Theorem \ref{teo kovalev maldonado}) controlled $H$-slope and round $H$-sections properties are strictly related:

\begin{theorem}\label{theorem mM se e solo se round}
Let $\varphi : \Hn\to \R$ be an $H$-convex function. The following conditions are equivalent:
\begin{itemize}
\item[a.] $\varphi$ is an $H$-differentiable function, with bounded $H$-sections and   controlled $H$-slope;

\item[b.] $\varphi$ has round $H$-sections.
\end{itemize}
Moreover, the constants $K_0$ and $K_1$ in (\ref{round section})  and in (\ref{H mM}) are related, and they depend only on  $\varphi$.
\end{theorem}
\psn
\noindent \textbf{Proof.} Let a. be true.
Let $S^H_\varphi(\xi_0,  s)$ be a bounded $H$-section, and let $R=\max\left\{d_g(\xi,\xi_0):\ \xi\in \overline{S^H_\varphi(\xi_0,  s)}\right\}$. Pick a point $\xi'$ such that $d_g(\xi',\xi_0)=R;$ then, $\xi'\in \partial {S^H_\varphi(\xi_0,  s)}$ and $\xi'=\xi_0\circ \exp v'$. From the $H$-convexity of $\varphi_{\xi_0}$ on $\H^n,$ we have that
$$
\varphi_{\xi_0}(\xi_0\circ\exp (v'/K_1))\le \left(1-\frac{1}{K_1}\right)\varphi_{\xi_0}(\xi_0)+\frac{1}{K_1}\varphi_{\xi_0}(\xi')=\frac{s}{K_1},
$$
where $K_1$ is as in (\ref{H mM}). Now,  for every $\xi\in H_{\xi_0}$ such that $d_g(\xi,\xi_0)=\frac{R}{K_1},$ by (\ref{H mM}) we have
$$
\varphi_{\xi_0}(\xi)\le M_\varphi^H\left(\xi_0,\frac{R}{K_1} \right)\le  K_1 m_\varphi^H\left(\xi_0,\frac{R}{K_1} \right)\le
\varphi_{\xi_0}(\xi')\le s.
$$
Hence,
$$
B\left(\xi_0,\frac{R}{K_1}\right)\subset S^H_\varphi(\xi_0,  s)\subset B\left(\xi_0,{R}\right).$$

Suppose now that condition b. holds true.
Proposition \ref{round implica Hdiffenziabilita} entails that $\varphi$ is $H$-differentiable. Consider $K_0$ as in \eqref{round section},   and fix $\xi\in\Hn$ and $r>0$: we have to prove (\ref{H mM}), where $K_1$ is uniform, i.e. it does not depend on $\xi$ and $r$. Set $s=m^H_\varphi(\xi,r)$ and define
$$
{\mathcal R}= \left\{R'>0:\  B_g(\xi, K_0R')\cap H_\xi\subset \overline{S^H_\varphi(\xi,s)}\subset B_g(\xi,R')\cap H_\xi\right\}.
$$
Since $\varphi$ has round $H$-sections, ${\mathcal R}$ is not empty. Set $R=\min {\mathcal R}$; trivially, $R=r,$ and
$$
\varphi(\xi\circ\exp(K_0v))-\varphi(\xi)-\nabla_H\varphi(\xi)\cdot (K_0v)\le s,
$$
for every $v\in V_1,\ \|v\|=R.$
The two relations above imply that
\begin{equation}\label{dim 2}
M_\varphi^H(\xi,K_0r)\le m_\varphi^H(\xi,r).
\end{equation}
Take $\alpha\in \N$ such that $K_0>2^{-\alpha},$ and note that relation (\ref{pro 1}) implies
$$
M_\varphi^H(\xi,R_1)\le C M_\varphi^H(\xi,R_1/2),
$$
for every $R_1>0,$ where $C$ depends only on $K_0$.
By iterating this relation, we obtain
$$
M_\varphi^H(\xi,r)\le C M_\varphi^H(\xi,2^{-1}r)\le C^2 M_\varphi^H(\xi,2^{-2}r)\le \ldots\le C^\alpha M_\varphi^H(\xi,2^{-\alpha}r)\le
 C^\alpha M_\varphi^H(\xi,K_0 r).$$
This last inequality, together with (\ref{dim 2}), leads to the assertion, with $K_1=C^{-\alpha}$  in (\ref{H mM}).
\hfill$\square$\psn

In the next result we investigate the properties of the function ${m}_\varphi^H,$ in order to shed some light on a finer behaviour of the $H$-sections.
\begin{proposition}\label{proposizione proprieta m H}
Let $\varphi :\Hn \to \R$ be an  H-differentiable and strictly H-convex function.
For every fixed $\xi\in\Hn$, the function  $r\mapsto  m_\varphi^H(\xi,r)$ is strictly increasing, continuous, and it goes to $+\infty,$ if $r\to +\infty$.
Then, the function $ m_\varphi^H(\xi,\cdot):[0,+\infty)\to[0,+\infty)$ is one-to-one and onto, and its inverse is defined on $[0,+\infty).$
A similar result holds for the function $ {M}^H_\varphi.$
\end{proposition}
\psn
\noindent \textbf{Proof.}
For every $\xi\in\Hn$, $r> 0$ and $v \in V_1,$ with $\|v\|=1$, set
$$
\widehat m_\varphi^H(\xi,v,r)=\min \{\varphi(\xi\circ \exp rv), \varphi(\xi\circ \exp (-r)v)\}.
%\min_{\{w\in V_1:\ w=\pm r v \}} \left(\varphi(\xi\circ \exp w)-\varphi(\xi)-\nabla_H \varphi(\xi)\cdot w\right).
$$
The function $\widehat m_\varphi^H$ is continuous , and strictly increasing w.r.t. $r,$ since $\varphi$ is strictly $H$-convex; thus,
$$
\widehat m_\varphi^H(\xi,v,r)<\widehat m_\varphi^H(\xi,v,r'), \qquad \forall \ 0\le r<r'.
$$
Hence, by the Berge Maximun Theorem (see, for instance, \cite{AliBor2002})
$m^H_\varphi$ is continuous, and
$$
{m}^H_\varphi(\xi,r)\le {m}^H_\varphi(\xi,r'),\qquad \forall \ 0\le r<r'.
$$
Let us show that the previous inequality is strict.
The set $\{v\in V_1:\  \|v\|=1\}$ is compact, and $\widehat m_\varphi^H(\xi, \cdot, \cdot)$ is continuous, then there exist $v$ and $v'$  such that
$\widehat m_\varphi^H(\xi,v,r)={m}^H_\varphi(\xi,r)$ and  $ \widehat m_\varphi^H(\xi,v',r')={m}^H_\varphi(\xi,r')$. This implies that
$$
{m}^H_\varphi(\xi,r)=
\widehat m_\varphi^H(\xi,v,r)\le  \widehat m_\varphi^H(\xi,v',r)<
 \widehat m_\varphi^H(\xi,v',r')={m}^H_\varphi(\xi,r').
$$
Let us show that $m^H_\varphi(\xi, \cdot)$ is unbounded, for every $\xi.$ Suppose, by contradiction, that there exists $L=L(\xi)>0$ such that ${m}^H_\varphi(\xi,r)\le L$ for every $r\ge 0.$ From the continuity of the function $v\mapsto \widehat m^H_\varphi(\xi,v,r),$  for every $r$ there exists $v_r$, with $\|v_r\|=1$,  such that $m^H_\varphi(\xi,r)=
\widehat m_\varphi^H(\xi,v_r,r).$ Let $r_n\to +\infty;$ then, there exists $\{v_{r_{n_k}}\}$ such that  $v_{r_{n_k}}\to \overline{v}.$ We have that
$$
\lim_{k\to +\infty}\widehat m_\varphi^H(\xi,v_{r_{n_k}},r_{n_k})=\lim_{k\to +\infty}\widehat m_\varphi^H(\xi,\overline{v},r_{n_k})=+\infty,
$$
contradicting the assumption that ${m}^H_\varphi(\xi,r)=\widehat m^H_\varphi(\xi,v_r,r)\le L$ for every $r>0.$
\hfill$\square$\psn

\section{Engulfing property for $H$-sections of $H$-convex functions}

This section is devoted to the study of the engulfing property $E(H,K)$ for the $H$-sections of an $H$-convex function. Our notion is different when compared with the one introduced by Capogna and Maldonado, and it generalizes the usual notion in the literature (see for example \cite{Gu2001}); however, we will see that these notions are equivalent (see Proposition \ref{propo soft equivalenti}). In the second part of the section we prove that a sufficient condition for a function to satisfies the engulfing property $E(H,K)$ is to have the  round $H$-sections property, or, equivalently, the controlled $H$-slope  (see Theorem \ref{theorem mM se e solo se round}). Finally, we will show, with an example, that the previous mentioned condition is only sufficient.

Let us start with our notion of engulfing for $H$-convex functions defined in $\H^n.$
\begin{definition}
Let  $\varphi :\Hn \to \R$ be an $H$-convex function.
We say that
$\varphi$ satisfies
the \emph{engulfing property} $E(H,K)$ (shortly, $\varphi \in E(H,K)$) if there exists $K >  1$
such that, for any $\xi\in\Hn$ and
$s > 0$, if $\xi'\in S_\varphi^H(\xi,p,s)$ with $p\in \partial_H \varphi(\xi)$, then
$$
S_\varphi^H(\xi,p,s)\cap H_{\xi'}\subset  S_\varphi^H(\xi',q,Ks)\cap H_{\xi},
$$
for every $q\in \partial_H \varphi(\xi')$.
\end{definition}
As a matter of fact, as mentioned previously, in \cite{CapMal2006} a slightly different definition of engulfing is investigated in the framework of Carnot groups; if $\textbf{G}=\H^n,$ it can be stated as follows:
\begin{gather}
\exists K>1: \ \texttt{\rm for every}\ \xi,\, \xi'\in\H^n\ \texttt{\rm and}\ s>0,\ \texttt{\rm if}\ \xi'\in S_\varphi^H(\xi,s),\ \texttt{\rm then}\ \xi\in  S_\varphi^H(\xi',Ks)\tag*{$\Diamond$}
\end{gather}
(we will refer to $\Diamond_K$ in case the constant $K$ plays a role).
Trivially, $\varphi\in  E(H,K)$ implies that $\varphi$ satisfies $\Diamond_K$.
The condition $\Diamond$ is essentially one-dimensional, as proved in the next
\begin{proposition}\label{prop Capogna Maldonado}(see \cite{CapMal2006}).
Let $\varphi:\Hn\to\R$ be a strictly $H$-convex and $H$-differentiable function. The function $\varphi$ satisfies $\Diamond_K$  if and only if for every $\xi\in\Hn$ and $v\in V_1$ the function $\varphi_{\xi,v}:\R\to\R$
 satisfies  condition ii. in Theorem \ref{th:softR}.
\end{proposition}

The following characterization provides an $\Hn$-version of the result in Theorem \ref{th:softR}:
\begin{proposition}\label{propo soft equivalenti}
Let $\varphi:\Hn\to\R$ be a strictly $H$-convex function. The following are equivalent:
\begin{itemize}
\item[i.] $\varphi$  satisfies the engulfing property $E(H,K),$ for some $K>1;$

\item[ii.]  $\varphi$ satisfies condition $\Diamond_{K'},$ for some $K'>1;$

\item[iii.] there exists a constant $K''> 1$ such that, for any $\xi\in\Hn,\ \xi'\in H_\xi,$ for any $p\in \partial_H\varphi(\xi)$ and $q\in \partial_H\varphi(\xi'),$
\begin{eqnarray*}
&&\frac{K''+1}{K''}\left(\varphi(\xi')-\varphi(\xi)-p\cdot(\texttt{\rm Pr}_1(\xi')-\texttt{\rm Pr}_1(\xi))\right)\\
&&\qquad\qquad \le(q-p)\cdot(\texttt{\rm Pr}_1(\xi')-\texttt{\rm Pr}_1(\xi))\\
&&\qquad\qquad\le(K''+1)\left(\varphi(\xi')-\varphi(\xi)-p\cdot(\texttt{\rm Pr}_1(\xi')-\texttt{\rm Pr}_1(\xi))\right).
\end{eqnarray*}
\end{itemize}
In particular, if any of the conditions above holds, $\varphi$ is $H$-differentiable.
\end{proposition}
\noindent \textbf{Proof.}
Trivially, i. implies ii., and one can take $K'=K.$
Let us show that ii. implies i.  Let $\xi'=\xi\circ \exp v$ be a point in $S^H_\varphi(\xi, p,s),$ and consider the convex function
 $\widehat\varphi_{\xi,v}:\R\to\R$ defined as in (\ref{def widehat varphi xi0 v}).
Note that
$$
S^H_\varphi(\xi,p,s)\cap H_{\xi'}=\{\xi\circ \exp sv:\, s\in S_{\widehat\varphi_{\xi,v}}(0,p\cdot v,s)\},
$$
and the function $\widehat\varphi_{\xi,v}$ satisfies condition ii. in Th. \ref{th:softR} with constant $K'.$ From Theorem 1 in \cite{CaPi2020}, $\widehat\varphi_{\xi,v}\in C^1(\R).$ Since this holds for every $\xi, v,$  from Remark \ref{rm:subdifferential} $\varphi$ is $H$-differentiable everywhere and $\partial_H\varphi(\xi)=\{\nabla_H\varphi(\xi)\}.$ Moreover, from Theorem 5.1 in \cite{CUFoMa2008}, the function $\widehat\varphi_{\xi,v}$ satisfies the engulfing condition with constant $2K'(K'+1).$  This is equivalent to say that
\begin{eqnarray}
&&\{\a\in \R:\, \widehat\varphi_{\xi,v}(\a)-\widehat\varphi_{\xi,v}(0)-\widehat\varphi_{\xi,v}'(0)\a <s\}\nonumber\\
&&\qquad\quad\subset \{\a\in \R:\, \widehat\varphi_{\xi,v}(\a)-\widehat\varphi_{\xi,v}(1)-\widehat\varphi_{\xi,v}'(1)(\a -1)<2K'(K'+1)s\}.\label{eq:inclusion}
\end{eqnarray}
From \eqref{eq:inclusion}, we get that
\begin{eqnarray*}
&&\{\a\in \R:\, \varphi(\xi\circ \exp \a v)-\varphi(\xi)-\nabla_H\varphi(\xi)\cdot v \a<s\}\\
&&\qquad\quad\subset \{\a\in \R:\, \varphi(\xi\circ \exp \a v)-\varphi(\xi\circ \exp v)-\nabla_H\varphi(\xi\circ \exp v)\cdot v(\a-1)<2K'(K'+1)s\},
\end{eqnarray*}
i.e.,
$\varphi$ is in $E(H,2K'(K'+1)).$

In order to prove that ii. implies iii., let $\xi'=\xi\circ \exp v$ and consider the convex function $\widehat\varphi_{\xi,v}.$ Note that
$p\cdot v\in \partial \widehat\varphi_{\xi,v}(0)$ and $q\cdot v\in \partial \widehat\varphi_{\xi,v}(1).$ Then, by applying Proposition 2.1 in \cite{CaPi2020}, we have that iii. holds with
$K''=K'.$ To conclude, let us show that iii. implies ii.  Take $\xi'=\xi\circ \exp v\in S^H_\varphi(\xi, p,s),$ where $p\in \partial_H\varphi(\xi),$ and let $q\in \partial_H\varphi(\xi').$ Then,
$$
\varphi(\xi)-\varphi(\xi\circ\exp v)-q\cdot (-v)\le \frac{K''}{K''+1}(q-p)\cdot v.
$$
The second inequality in iii. gives
$$
(p-q)\cdot (-v)\le (K''+1)(\varphi (\xi\circ \exp v)-\varphi(\xi)-p\cdot v) \le (K''+1)s.
$$
Then,
$$
\varphi(\xi)-\varphi(x\circ \exp v)-q\cdot (-v)\le K'' s,
$$
thus, $\xi\in S^H_\varphi(\xi\circ \exp v,q,K''s),$ i.e., condition $\Diamond_{K''}$ is fulfilled.
\hfill$\square$\psn

Let us recall that a set-valued map $T:\H^n\to \mathcal{P}(V_1)$ is said to be $H$-monotone if, for all $\xi\in \H^n,$  $\xi'\in  H_\xi,$ $p\in T(\xi),$ $q\in T(\xi'),$ then
$$
(q-p)\cdot ( \texttt{\rm Pr}_1(\xi')-\texttt{\rm Pr}_1(\xi))\ge 0
$$
(here $V_1\cong \R^n$). In particular, if $\varphi$ is an $H$-convex function, then the $H$-subdifferential map $\partial_H\varphi$ is an $H$-monotone set-valued map (see \cite{CaPi2014}). The property iii. above requires, in fact, a stronger control on the $H$-monotonicity, both from below and from above.

Let us now state the following crucial result, that provides a sufficient condition for  $E(H,K)$ via the round $H$-sections property; the relationship between round $H$-sections, or, equivalently, controlled $H$-slope, and the engulfing property corresponds to the similar one in $\R^n$, for $n\ge 2$:
\begin{theorem}\label{round implica engulfing}
If $\varphi:\Hn\to\R$ is an $H$-convex function with round $H$-sections, then $\varphi$ satisfies  the  engulfing property $E(H,K),$ where $K$ depends only on $K_0$ in  (\ref{round section}).
\end{theorem}
\psn
\noindent \textbf{Proof.}
Since $\varphi$ has round $H$-sections, Proposition \ref{round implica Hdiffenziabilita} implies that $\varphi$ is  strictly $H$-convex and $H$-differentiable. Let
 $\xi'\in S_\varphi^H(\xi,s)$ be such that $\xi'=\xi\circ\exp (r'v)$ for some $v$ in $V_1,$ with $\|v\|=1$ and $r'>0$; we will prove that $\xi\in  S_\varphi^H(\xi',Ks)$ where $K$ depends only on $K_0$ in (\ref{round section}).

Let $R$ be such that
 \begin{equation}\label{round section 1}
B_g(\xi, K_0R)\cap H_\xi\subset S^H_\varphi(\xi,s)\subset B_g(\xi,R)\cap H_\xi.
\end{equation}
Since  $S_\varphi^H(\xi,s)$ is bounded, let us consider
$$
r^\partial=\max\left\{r\ge 0:\ \xi\circ\exp (rv)\in \overline{S_\varphi^H(\xi,s)}\right\},\qquad \xi^\partial=\xi\circ\exp (r^\partial v)\in \partial S_\varphi^H(\xi,s).
$$
Hence,
 \begin{equation}\label{dim 11}
K_0R\le r^\partial\le R,
\end{equation}
and $0<r'\le r^\partial$. From the $H$-monotonicity of the map $\xi\mapsto \partial_H\varphi(\xi)$  we have that

\begin{equation}\label{dim 11 7}
0\le\left(\nabla_H\varphi(\xi')-\nabla_H\varphi(\xi)\right)\cdot v\le \left(\nabla_H\varphi(\xi^\partial)-\nabla_H\varphi(\xi)\right)\cdot v.
\end{equation}
Let us introduce the function $\Phi:\R\to\R$ defined by
$$
 \Phi(\a)=\widehat\varphi_{\xi, v}(\alpha),\qquad \forall \a\in\R;
$$
this function is strictly convex,  with  $\Phi(0)=\Phi'(0)=0$. Let us consider the function $\Pi:\R\to\R$ defined by
$$
 \Pi(\a)=\Phi(r^\partial)+\Phi'(r^\partial)(\a-r^\partial),\qquad \forall\a\in\R;
$$
clearly, it represents the tangent to the graph of $\Phi$ at $(r^\partial, \Phi(r^\partial))$ with $\Phi(r^\partial)>0$ and $\Phi'(r^\partial)>0;$ hence we have
\begin{eqnarray}
 \Pi(\a)&=&\varphi(\xi\circ\exp(r^\partial v))-\varphi(\xi)-\nabla_H\varphi(\xi)\cdot vr^\partial+\nonumber\\
 &&\quad+\left(\nabla_H\varphi(\xi\circ\exp(r^\partial v))-\nabla_H\varphi(\xi)\right)\cdot v(\a-r^\partial )\nonumber\\
 &=&\varphi(\xi\circ\exp(\a v))-\varphi(\xi)-\nabla_H\varphi(\xi)\cdot v\a +\nonumber\\
 &&\quad-\left(\varphi(\xi\circ\exp(\a v))-\varphi(\xi\circ\exp(r^\partial v))-\nabla_H\varphi(\xi\circ\exp(r^\partial v))\cdot v(\a-r^\partial )\right)\nonumber
\end{eqnarray}
Since $\varphi$ is $H$-convex, the previous equalities and (\ref{def varphi xi}) give
\begin{eqnarray}
 \Pi(\a)&\ge&\left(\nabla_H\varphi(\xi^\partial )-\nabla_H\varphi(\xi)\right)\cdot v(\a-r^\partial ),\quad\forall \a\label{dim 12}\\
 \Pi(\a)&\le&\varphi_\xi(\xi\circ\exp(\a v))\label{dim 13},\quad\forall \a
\end{eqnarray}
From (\ref{dim 11 7}) and (\ref{dim 12}) we get
%\begin{equation}\label{dim 14}
$$\left(\nabla_H\varphi(\xi')-\nabla_H\varphi(\xi)\right)\cdot v
\le
\left(\nabla_H\varphi(\xi^\partial)-\nabla_H\varphi(\xi)\right)\cdot v
\le \frac{\Pi(2r^\partial)}{r^\partial}.
%\end{equation}
$$
The inequality above, together with (\ref{dim 11}) and  (\ref{dim 13}), give
$$
\left(\nabla_H\varphi(\xi')-\nabla_H\varphi(\xi)\right)\cdot v
\le
\frac{\varphi_\xi(\xi\circ\exp(2r^\partial v))}{K_0R}.
$$
The $H$-convexity of $\varphi$ and (\ref{dim 11}) imply that $\varphi_\xi(\xi\circ\exp(2r^\partial v))\le \varphi_\xi(\xi\circ\exp(2R v))$; hence
 we obtain
\begin{equation}\label{dim 15 bis}
\left(\nabla_H\varphi(\xi')-\nabla_H\varphi(\xi)\right)\cdot v
\le
\frac{\varphi_\xi(\xi\circ\exp(2R v))}{K_0R}.
\end{equation}
Let us consider $\alpha\in \N$ such that $K_0>2^{-\alpha}.$ By iterating  relation (\ref{pro 1}), we obtain
$$
\varphi_\xi(\xi\circ\exp(2R v))\le C \varphi_\xi(\xi\circ\exp(R v))\le C^{1+\alpha} \varphi_\xi(\xi\circ\exp(R2^{-\alpha} v))
\le C^{1+\alpha} \varphi_\xi(\xi\circ\exp(K_0R v)),$$
where $C$ depends only on $K_0$. The previous inequality,  and relations (\ref{round section 1}) and (\ref{dim 15 bis}), give
\begin{equation}\label{dim 16}
\left(\nabla_H\varphi(\xi')-\nabla_H\varphi(\xi)\right)\cdot v
\le C^{1+\alpha} \frac{s}{K_0R}.\end{equation}
At this point, since  $S_\varphi^H(\xi,s)$ is open, there exists $\tilde\xi=\xi\circ\exp(\tilde r v)\in S_\varphi^H(\xi,s)$ with $\tilde r<0$.  Taking into account that $\xi'\in S_\varphi^H(\xi,s)$ and  $\tilde \xi\in S_\varphi^H(\xi,s)$, and using (\ref{dim 16}) and (\ref{round section 1}), we obtain
\begin{eqnarray*}
&&\varphi(\tilde\xi)-\varphi(\xi')-\nabla_H\varphi(\xi')\cdot ({\rm Pr}_1(\tilde\xi)-{\rm Pr}_1(\xi'))=\\
&&\qquad= \varphi(\tilde\xi)-\varphi(\xi)-\nabla_H\varphi(\xi)\cdot ({\rm Pr}_1(\tilde\xi)-{\rm Pr}_1(\xi))+\\
&&\qquad\qquad-\left(\varphi(\xi')-\varphi(\xi)-\nabla_H\varphi(\xi)\cdot ({\rm Pr}_1(\xi')-{\rm Pr}_1(\xi))\right)+\\
&&\qquad\qquad-\left(\nabla_H\varphi(\xi')-\nabla_H\varphi(\xi)\right)\cdot({\rm Pr}_1(\tilde \xi)-{\rm Pr}_1(\xi'))\\
& &\qquad < s+\left(\nabla_H\varphi(\xi' )-\nabla_H\varphi(\xi)\right)\cdot v(r'-\tilde r) \\
& &\qquad \le s\left(1+ \frac{2C^{1+\alpha }}{K_0}\right).
\end{eqnarray*}
This implies that $\tilde\xi\in S_\varphi^H(\xi',Ks)$, with $K=1+ \frac{2C^{1+\alpha }}{K_0}$: since $\xi$ belongs to the horizontal segment which joins $\tilde \xi$ and $\xi',$ and since $S_\varphi^H(\xi',Ks)$ is $H$-convex, then  $\xi\in S_\varphi^H(\xi',Ks).$ By Proposition \ref{propo soft equivalenti} the assertion is proved.
\hfill$\square$\psn

The following example  is crucial in order to shed some light on the relationship between round sections and engulfing; indeed, it shows that the converse of the previous theorem fails. The idea is taken from an example due to Wang (see \cite{Wa1995}) and set in $\R^2;$ we adapt his idea to the case of the first Heisenberg group $\H$.
\begin{example}\label{example Wang}{\rm
Consider the following differentiable and strictly convex function $u:\R^2\to \R,$
\begin{equation}\label{ex Wang}
u(x, y) =
\left\{
\begin{array}{ll}
\displaystyle x^4+\frac{3 y^2}{2x^2} &|y|\le |x|^3\\
\displaystyle \frac{1}{2}x^2|y|^{2/3}+2|y|^{4/3}& |y|>|x|^3.\\
\end{array}\right. \end{equation}
The Monge-Amp\`{e}re measure $\mu_u$ (we recall that $\mu_u$ is defined by $\mu_u(E)=|\partial u(E)|$ for every Borel set $E\subset\R^2)$ is absolutely continuous with respect to the
Lebesgue measure $|\cdot|$, and it verifies the condition $\mu_\infty$, i.e. for any $\delta_1\in (0,1)$ there exists $\delta_2\in (0,1)$ such that:
 for every section $S_u(z,s)$, with $z\in\R^2$, and for every Borel set $B\subset S_u(z,s),$
 $$
 \frac{|B|}{|S_u(z,s)|}<\delta_2\quad\Rightarrow\quad
 \frac{\mu_u(B)}{\mu_u(S_u(z,s))}<\delta_1
 $$
(see Definition 3.7 in \cite{KoMa2005}).
This condition $\mu_\infty$ is stronger than the (DC)-doubling property (see, for example, relation (3.1.1) in \cite{Gu2001}), i.e.,  there exist constants $\alpha\in (0,1)$ and $C>1$ such that
$$
\mu_{u}(S_{u}(z,s))\le C\mu_{u}(\alpha S_{u}(z,s)),
$$
for every $z,s>0$ (here $\alpha S_{u}(z,s)$ denotes the open convex set obtained by $\alpha$-contraction of $S_{u}(z,\t)$ with respect to its center of mass).
In \cite{GuHu2000} and \cite{FoMa2002} it was shown that the (DC)-doubling property of the measure $\mu_{u}$ is equivalent to the engulfing property of the function $u.$ Therefore, $u$ satisfies the engulfing property.

Since the second derivative of $u$ w.r.t. $x_2$
is unbounded near the origin, so is $\|D^2 u\|$; thus, $u$ is not quasiuniformly
convex (see i. in Theorem \ref{teo kovalev maldonado} and \cite{KoMa2005} for further details).
However, a simpler argument can be advanced to prove that $u$ is not quasiuniformly
convex, that is one can show that
$u$ has not controlled slope in (\ref{def controlled slope}): in order to do that,
we only remark that, taking into account that $u(0,0)=0$ and $\nabla u(0,0)=(0,0)$, we have, for large $r,$
 \begin{eqnarray*}
m_u((0,0),r)&=&\min_{\{z\in\R^2:\ \|z\|=r\}} u(z)\le \varphi(0,r)=2r^{4/3}\\
M_u((0,0),r)&=&\max_{\{z\in\R^2:\ \|z\|=r\}} u(z)\ge \varphi(r,0)=r^4.
 \end{eqnarray*}
Now let us consider the function $\varphi:\H\to\R$ defined by
$\varphi(x,y,t)=u(x,y),$ for all $(x,y,t)\in \H.$ This function $\varphi$ is $\R^3$-convex, and hence $H$-convex. Since
$$
m_u((0,0),r)=m_\varphi^H((0,0,0),r),\qquad M_u((0,0),r)=M_\varphi^H((0,0,0),r).
$$
$\varphi$ has not controlled $H$-slope, and hence has not round $H$-sections.
However, since
$$
(x,y)\in S_{u}((x_0,y_0),s)\qquad \Longleftrightarrow \qquad (x,y,t)\in S_{\varphi}^H((x_0,y_0,t_0),s),
$$
it is easy to see that $\varphi$ enjoys the engulfing property.}

\end{example}

\section{$\Hn$-sections of $H$-convex functions and their engulfing properties}

In this section we will present our new definition of section in $\H^n.$ First of all, we will prove that these $\Hn$-sections have topological dimension $2n+1$, thereby allowing to construct a topology in $\H^n,$ as we will see in the next Section 6. In the second part, we introduce the condition of engulfing $E(\Hn,K)$ for  these new $\Hn$-sections. It will not be a surprise that $\varphi\in E(\Hn,K)$ implies that
 $\varphi\in E(H,K)$, while the converse implication is very hard and mysterious (at least to us). In order to shed some light on this, let us focus our attention on the functions having round $H$-sections, or, equivalently, controlled $H$-slope. As we will see, some technical estimates allow us to prove the first part of our main result in Theorem \ref{intro Teo}.

 Let us start with our new notion of $\Hn$-section:

\begin{definition}\label{definizione section Heis}
 Let $\varphi:\Hn \to\R$ be an H-convex function and let us fix $\xi_0\in\Hn.$ For a
given $s>0,$ an \emph{$\Hn$-section} of $\varphi$ at height $s$, with $p_0\in \partial_H \varphi(\xi_0)$, is the set
\begin{equation}\label{def section Heis}
\S_\varphi^\Hn(\xi_0,p_0,s)=\bigcup_{
\begin{array}{l}
\xi_1\in S_\varphi^H(\xi_0,p_0,s),\ p_1\in \partial_H \varphi(\xi_1),\\
\xi_2\in S_\varphi^H(\xi_1,p_1,s),\ p_2\in \partial_H \varphi(\xi_2)\\
\end{array}
} \hskip -2.1truecm S_\varphi^H(\xi_2,p_2,s).
\end{equation}
In case $\varphi$ is $H$-differentiable at $\xi_0$, we will denote the $\Hn$-section at $\xi_0$ with height $s$ by $\S_\varphi^\Hn(\xi_0,s)$, for
short.
\end{definition}
Let us spend a few words on the definition above. Lemma 1.40 in the fundamental book by Folland and Stein \cite{FoSt1982} guarantees that, in every stratified group $(\textbf{G},\circ)$ with homogeneous norm $\|\cdot\|_\textbf{G}$, there exists a constant $C>0$ and an integer $k\in\N$ such that any $\xi\in\textbf{G}$ can be expressed as $\xi=\xi_1\circ\xi_2\circ\ldots\circ\xi_k$, with $\xi_i\in \exp(V_1)$ and $\|\xi_i\|_\textbf{G}\le C \|\xi\|_\textbf{G},$ for every $i$. If $\textbf{G}=\H^n,$ the mentioned $k$ is exactly 3, for every $n\ge 1.$ In other words,  every point $\xi\in\Hn$ can be reached from the origin $e$ following a path of three consecutive horizontal segments. The idea behind Definition \ref{definizione section Heis} takes inspiration from this result, in view of providing a family of sets with nonempty interior. Let us define, for every $\xi\in\Hn$ and $r>0$,
\begin{equation}\label{balls tilde}
\widetilde{B}(\xi,r)=\left\{\xi'\in\Hn:\ \xi'=\xi\circ\exp(v_1)\circ\exp(v_2)\circ\exp(v_3);\ v_i\in V_1,\ \|v_i\|\le r\right\}.
\end{equation}
Clearly, $\delta_\lambda\left(\widetilde{B}(e,r)\right)=\widetilde{B}(e,\lambda r),$ and the associated distance $\widetilde d$ in $\Hn$ is left-invariant and homogeneous; hence, it is bi-Lipschitz equivalent to $d_g$ and to any other left-invariant and homogeneous distance in $\Hn$.
 Moreover, due to the  Folland--Stein Lemma, we have that, for every $\xi\in \Hn$ and $r>0$,
\begin{equation}\label{inclusione palle}
\widetilde{B}(\xi,r)\subset {B}_g(\xi,3r)\subset \widetilde{B}(\xi,3Cr),
\end{equation}
where $C$ is the constant in the mentioned lemma.

Let us prove the first fundamental property of the $\Hn$-sections, i.e. that $\S_\varphi^\Hn(\xi_0,p_0,s)$ has a topological dimension equal to $2n+1$.
\begin{proposition}
Let  $\varphi :\Hn \to \R$ be an $H$-convex function. Then, for every $\xi_0\in\Hn,\ p_0\in\partial_H \varphi(\xi_0)$ and $s>0$, there exists $r>0$ such that
$$
B_g(\xi_0,r)\subset \S_\varphi^\Hn(\xi_0,p_0,s).$$
\end{proposition}
\psn
\noindent \textbf{Proof.} Without loss of generality, we set $\xi_0=e$. Let $r>0,$ and assume that $B_g(\xi_0,r)\not\subset \S_\varphi^\Hn(\xi_0,p_0,s)$. Denote by $K$ the compact set
$$
K=\overline{B_g(\xi_0,r)\cap \S_\varphi^\Hn(\xi_0,p_0,s).}
$$
Since the $H$-subdifferential map $\partial_H\varphi$ brings compact sets into compact sets (see, for instance, Proposition 2.1 in \cite{BaCaKr2015}), there exists a constant $R=R(s,r)$ such that
\begin{equation}\label{prop new 1}
\partial_H \varphi(K)\subset B^{\R^{2n}}(0,R).
\end{equation}
Moreover, since $\varphi$ is locally Lipschitz (see Theorem 1.2 in \cite{BaRi2003}), there exists a constant $L=L(s,r)$ such that
\begin{equation}\label{prop new 2}
|\varphi(\xi)-\varphi(\xi')|\le L d_g(\xi,\xi'),\qquad \forall \xi,\xi'\in K.
\end{equation}
Define $r=\min\left(\frac{s}{(L+2R)C},r\right),$ where $C$ is the constant in the Folland--Stein Lemma; we will prove that $B_g(e,r)\subset \S_\varphi^\Hn(e,p_0,s).$ Take any $\xi\in  B_g(e,r)$; then,
 $\xi=\exp(v_1)\circ\exp(v_2)\circ\exp(v_3)$  for suitable $\{v_i\}_{i=1}^3\subset  V_1$ such that $\|v_i\|\le C \|\xi\|_g\le r,$ for every $i$. Set $\xi_i=\exp(v_1)\circ\ldots\circ\exp(v_i),$ $i=1,2,3.$ Then, by (\ref{prop new 1}) and (\ref{prop new 2}), we have
\begin{equation}\label{prop new 3}
\varphi(\xi_1)-\varphi(e)-p_0\cdot v_1 \le (L+R)\|v_1\|< s.
 \end{equation}
 Similarly, for every $i=2,3$ and $p_i\in\partial_H\varphi(\xi_i)$ we have
 \begin{equation}\label{prop new 4}
\varphi(\xi_i)-\varphi(\xi_{i-1})-p_i\cdot v_i \le (L+R)\|v_i\|< s.
 \end{equation}
From (\ref{prop new 3}) and  (\ref{prop new 4}) we get the claim.
\hfill$\square$

Starting from these $\H^n$-sections, we introduce the following engulfing property:
\begin{definition}
Let  $\varphi :\Hn \to \R$ be an $H$-convex function.
We say that
$\varphi$ satisfies
the \emph{engulfing property} $E(\Hn,K)$ if there exists $K > 1$
such that, for any $\xi\in\Hn,\ p\in \partial_H \varphi(\xi)$ and
$s > 0$, if $\xi'\in \S_\varphi^\Hn(\xi,p,s)$, then $$\S_\varphi^\Hn(\xi,p,s)\subset  \S_\varphi^\Hn(\xi',q,Ks),$$ for every $q\in \partial_H \varphi(\xi')$.
\end{definition}
This engulfing property $E(H,K)$ is related to this the engulfing property $E(\H^n,K)$ as well as condition $\Diamond$ is related to the following condition:
\begin{gather}\tag*{$\Diamondblack$}
\exists K>1: \ \texttt{\rm for every}\ \xi\in\Hn,\ p\in \partial_H \varphi(\xi) \
 \texttt{\rm and}\ s>0 \\
 \texttt{\rm if}\ \xi'\in \S_\varphi^\Hn(\xi,p,s),\ \texttt{\rm then}\ \xi\in  \S_\varphi^\Hn(\xi',q,K's) \ \texttt{\rm for every}\ q\in \partial_H \varphi(\xi').\nonumber
\end{gather}
We will refer to $\Diamondblack_K$ in case we need to specify the constant $K$ in the previous condition.

It is clear that
\begin{remark}\label{remark ciccio}
  If $\varphi$ satisfies the {engulfing property} $E(\Hn,K)$, then condition  $\Diamondblack_K$ holds.
\end{remark}
The converse of the previous remark is a delicate question: the aim of this section is, essentially, to prove that, under further conditions on $\varphi$, the converse of Remark \ref{remark ciccio} holds.

The relationship between conditions $\Diamond$ and $\Diamondblack$ is the following:
\begin{proposition}\label{proposizione se e solo se soft engulfing}
Let $\varphi :\Hn \to \R$ be an  H-convex function. Then $\varphi$ satisfies condition $\Diamond_K$
if and only if $\varphi$   satisfies condition $\Diamondblack_K$.
\end{proposition}
\psn
\noindent \textbf{Proof.} If $\varphi$   satisfies  $\Diamondblack_K$, it is clear that $\Diamond_K$ holds. Let us prove the converse. Take any $\xi'\in\S_\varphi^\Hn(\xi,p,s)$, i.e.
$\xi'=\xi\circ \exp(v_1)\circ \exp(v_2)\circ\exp(v_3),$ with $v_i\in V_1$ and with
\begin{eqnarray*}
&&\xi_1:=\xi\circ \exp(v_1)\in S_\varphi^H(\xi,p,s),\\
&&\xi_2:=\xi\circ \exp(v_1)\circ \exp(v_2)\in S_\varphi^H(\xi_1,p_1,s),\qquad \texttt{\rm with}\ p_1\in \partial_H \varphi(\xi_1),\\
&&\xi'\in S_\varphi^H(\xi_2,p_2,s),\qquad \qquad\qquad\qquad\qquad\qquad\quad\texttt{\rm with}\ p_2\in \partial_H \varphi(\xi_2);
\end{eqnarray*}
we have to show that $\xi\in  \S_\varphi^\Hn(\xi',q,K's)$, for every $q\in \partial_H \varphi(\xi')$. The assumption implies
\begin{eqnarray*}
&&\xi\in  S_\varphi^H(\xi_1,p_1,K's),\qquad \forall p_1\in \partial_H \varphi(\xi_1),\\
 && \xi_1\in S_\varphi^H(\xi_2,p_2,K's),\qquad \forall p_2\in \partial_H \varphi(\xi_2),\\
  && \xi_2\in S_\varphi^H(\xi',q,K's),\qquad\quad \forall q\in \partial_H \varphi(\xi').
\end{eqnarray*}
Hence, for every   $q\in \partial_H \varphi(\xi')$,
$$
\xi\in
\bigcup_{
\begin{array}{l}
\xi_2\in S_\varphi^H(\xi',q,K's),\ p_2\in \partial_H \varphi(\xi_2),\\
 \xi_1\in S_\varphi^H(\xi_2,p_2,K's),\ p_1\in \partial_H \varphi(\xi_1)
 \end{array}
} \hskip -1.7truecm
 S_\varphi^H(\xi_1,p_1, K's)=\S_\varphi^\Hn(\xi',q,K's).
$$
\hfill$\square$\psn
\medskip

Clearly, if $\varphi$ is a strictly $H$-convex function  satisfying the engulfing property $E(\Hn,K)$, then Remark \ref{remark ciccio},
Proposition \ref{proposizione se e solo se soft engulfing} and Proposition \ref{propo soft equivalenti} imply that
$\varphi$ is $H$-differentiable.

The next result will be crucial to our purposes:

\begin{proposition}\label{proposizione due inclusioni}
Let $\varphi :\Hn \to \R$ be an  H-differentiable and strictly H-convex function.
 Then,
 for every $ \xi\in\Hn,\ r>0$, we have
 \begin{eqnarray}
\bigcup_{
\begin{array}{l}
\xi_1\in S^H_\varphi(\xi,  m_\varphi^H(\xi,r))\\
 \xi_2\in S^H_\varphi(\xi_1,  m_\varphi^H(\xi_1,r))
 \end{array}
}\hskip -1.5truecm
S^H_\varphi(\xi_2,   m_\varphi^H(\xi_2 ,r))
&\subset&
\widetilde{B}(\xi,r)\label{gen incl 1}\\
&\subset&
 \bigcup_{
\begin{array}{l}
\xi_1\in S^H_\varphi(\xi,  M_\varphi^H(\xi,r))\\
 \xi_2\in S^H_\varphi(\xi_1,  M_\varphi^H(\xi_1,r))
 \end{array}
} \hskip -1.5truecm
 S^H_\varphi(\xi_2,   M_\varphi^H(\xi_2,r)).\label{gen incl 2}
 \end{eqnarray}

\end{proposition}
Let us emphasize that, despite its appearance, the first set in (\ref{gen incl 1}) is not an $\Hn$-section, since $m_\varphi^H(\xi_1,r)$, for  $\xi_1\in S^H_\varphi(\xi,  m_\varphi^H(\xi,r))$, and $m_\varphi^H(\xi_2,r)$, for  $\xi_2\in S^H_\varphi(\xi_1,  m_\varphi^H(\xi_1,r))$, are not fixed values. A similar comment holds for the set in (\ref{gen incl 2}).

\noindent \textbf{Proof of Proposition \ref{proposizione due inclusioni}.}
By the  inclusions in (\ref{inclusione con palla 2}), we easily have
 \begin{eqnarray*}
\bigcup_{
\begin{array}{l}
\xi_1\in S^H_\varphi(\xi,  m_\varphi^H(\xi,r))\\
 \xi_2\in S^H_\varphi(\xi_1,  m_\varphi^H(\xi_1,r))
 \end{array}
} \hskip -1.7truecm
S^H_\varphi(\xi_2,   m_\varphi^H(\xi_2,r))
&\subset&
\bigcup_{
\begin{array}{l}
\xi_1\in S^H_\varphi(\xi,  m_\varphi^H(\xi,r))\\
 \xi_2\in S^H_\varphi(\xi_1,  m_\varphi^H(\xi_1,r))
 \end{array}
} \hskip -1.7truecm
 B_g(\xi_2,r)\cap H_{\xi_2}\\
 &\subset&
\bigcup_{
\begin{array}{l}
\xi_1\in B_g(\xi,r)\cap H_{\xi}\\
 \xi_2\in B_g(\xi_1,r)\cap H_{\xi_1}
 \end{array}
} \hskip -1.5truecm
 B_g(\xi_2,r)\cap H_{\xi_2}\\
 &\subset&
 \bigcup_{
\begin{array}{l}
\xi_1\in B_g(\xi,r)\cap H_{\xi}\\
 \xi_2\in B_g(\xi_1,r)\cap H_{\xi_1}
 \end{array}
} \hskip -1.5truecm
 S^H_\varphi(\xi_2,   M_\varphi^H(\xi_2,r)) \\
 &\subset&
 \bigcup_{
\begin{array}{l}
\xi_1\in S^H_\varphi(\xi,  M_\varphi^H(\xi,r))\\
 \xi_2\in S^H_\varphi(\xi_1,  M_\varphi^H(\xi_1,r))
 \end{array}
} \hskip -1.7truecm
 S^H_\varphi(\xi_2,   M_\varphi^H(\xi_2,r))
\end{eqnarray*}
for every $ \xi\in\Hn,\ r>0$. Hence the assertion holds.
\hfill$\square$\psn
\medskip

In order to prove our main result concerning the engulfing property of the $\H^n$-sections, an extension to the Heisenberg case of the inequalities (\ref{B0}), (\ref{B4}) and (\ref{B2}) turns out to be quite useful:

\begin{proposition}\label{proposizione M e m in H zero}
Let $\varphi$ be a strictly $H$-convex function in $E(H,K)$. Then, for every $r\ge 0$ and $\xi\in\Hn,$ we have
\begin{eqnarray}
&&{M}_\varphi^H(\xi,2r)\le B_1
 {M}_\varphi^H(\xi,r),\label{B_1 in H zero}\\
 &&{m}_\varphi^H(\xi,2r)\le B_2
 {m}_\varphi^H(\xi,r),\label{B_2 in H zero}\\
 &&B_4 m_\varphi^H(\xi,r)\le m_\varphi^H(\xi,2r),\label{B_4 in H zero}
\end{eqnarray}
where $  B_1,\ B_2$ and $B_4$ depend only on $K,$ and $B_i>1$.
\end{proposition}
\psn
\noindent \textbf{Proof.} Proposition \ref{propo soft equivalenti} implies that $\varphi$ is $H$-differentiable and, if we consider its restriction to any horizontal segment, we obtain a strictly convex and differentiable function. To be precise,
 for every $\xi\in\Hn$ and $v\in V_1$ with $\|v\|=1$ the function $\widehat\varphi_{\xi,v}:\R\to\R$, defined as in (\ref{def widehat varphi xi0 v}), satisfies condition ii. in Th. \ref{th:softR}. By (\ref{B0}) in
 Remark \ref{remark} we obtain
$$
{M}_{\widehat\varphi_{\xi,v}}(0,2r)\le  B_1
{M}_{\widehat\varphi_{\xi,v}}(0,r),
$$
where $B_1$ depends only on $K$.
Hence we have
\begin{eqnarray*}
&&\max_{\{w\in V_1:\ w=\pm 2rv\}}\left(\varphi(\xi\circ\exp w)-\varphi(\xi)-\nabla_H \varphi(\xi)\cdot w\right)\le\\
  &&\qquad\qquad\le B_1
\max_{\{w\in V_1:\ w=\pm rv\}}\left(\varphi(\xi\circ\exp w)-\varphi(\xi)-\nabla_H \varphi(\xi)\cdot w\right);
\end{eqnarray*}
 taking the  maximum w.r.t. to  $v$, with $\|v\|=1,$ we obtain (\ref{B_1 in H zero}).

 A similar proof, via inequality (\ref{B4}) in Remark \ref{remark} and inequality (\ref{B2}) in Proposition \ref{prop new}, shows
 (\ref{B_2 in H zero}) and (\ref{B_4 in H zero}), respectively.
\hfill$\square$\psn

In the final part of this section we will prove our main result concerning the relationship between round $H$-sections and the engulfing property of the $\H^n$-sections. The proof will be quite technical, deserving a few previous estimates.

Let $\varphi:\H^n\to \R$ be an $H$-convex function with round $H$-sections (with constant $K_0$). Then,  $\varphi\in E(H,K)$, and has controlled $H$-slope (with constant $K_1$), where both $K,K_1$ depend on $K_0.$
Denote by $\gamma$ any positive integer such that
\begin{equation}\label{gamma}
 K_1\le B_4^\gamma.
\end{equation}
Thus, from
(\ref{H mM}),  and by iterating inequality (\ref{B_4 in H zero}), we obtain
$$
 M_\varphi^H(\xi,r)\le K_1  m_\varphi^H(\xi,r)\le \frac{K_1}{B_4^\gamma}m_\varphi^H(\xi,2^{\gamma}r)\le m_\varphi^H(\xi,2^{\gamma}r).
$$
Then, we have that
\begin{equation}\label{remark_gamma}
M_\varphi^H(\xi,r)\le m_\varphi^H(\xi,2^{\gamma}r),
\end{equation}
for every $r>0$ and $\xi\in\Hn$, where $\gamma>1$ depends only on $K_0$ in (\ref{round section}).

The next proposition holds:
\begin{proposition}\label{prop tre in R}
Let $\varphi :\H^n \to \R$ be a function with round $H$-sections (with $K_0$ as in (\ref{round section})).
Then, there exists a constant $C_1>0$ such that, if $\xi'\in S_\varphi^H (\xi,s)$, then
%\begin{equation}\label{C_1 in H}
$$
 m^H_\varphi(\xi',r)\le C_1 m^H_\varphi(\xi,r),
%\end{equation}
$$
for $r$ such that  $s= m^H_\varphi(\xi,r)$. The constant $C_1$ depends only on $K_0$.
\end{proposition}
\noindent \textbf{Proof.} Since $\varphi$ has round $H$-sections, it is
 strictly $H$-convex, $H$-differentiable and it satisfies
the engulfing property $E(H,K)$, where $K$ depends only on $K_0$.
Let $\xi'=\xi\circ\exp v\in S_\varphi^H (\xi,s),$ and set $r$ such that $s= m^H_\varphi(\xi,r)$:
clearly,
$$
\xi'\in S_\varphi^H (\xi, m^H_\varphi(\xi,r))\subset  S_\varphi^H (\xi, M^H_\varphi(\xi,r))\subset S_\varphi^H (\xi, {M}^H_\varphi(\xi,2r)).
$$
Moreover, since $\xi'\in S_\varphi^H (\xi, m^H_\varphi(\xi,r))$, by (\ref{inclusione con palla 2}) we have that
$$\xi'\circ\exp(\pm rv/\|v\|)\in B_g(\xi,2r)\cap H_\xi\subset  S_\varphi^H (\xi,{M}^H_\varphi(\xi,2r)).$$
Furthermore, since $\varphi\in E(H,K)$, we have that
$\xi'\in S_\varphi^H (\xi,{M}^H_\varphi(\xi,2r))$ gives
$$S_\varphi^H (\xi,{M}^H_\varphi(\xi,2r))\cap H_{\xi'}\subset  S_\varphi^H (\xi',K' {M}^H_\varphi(\xi,2r))\cap H_{\xi}.$$
This implies $\xi'\circ\exp(\pm rv/\|v\|)\in S_\varphi^H (\xi',K{M}^H_\varphi(\xi,2r))$. Now, by (\ref{B_1 in H zero}) and (\ref{round section}), we have
$$
\xi'\circ\exp(\pm rv/\|v\|)\in S_\varphi^H (\xi',K M^H_\varphi(\xi,2r))\subset
S_\varphi^H (\xi',K B_1 M^H_\varphi(\xi,r))\subset S_\varphi^H (\xi',K B_1 K_0  m^H_\varphi(\xi,r)),
$$
where $B_1$ depends only on $K_0$. Then,
$$
m_\varphi^H(\xi',r)\le
\varphi(\xi'\circ\exp(\pm rv/\|v\|))-\varphi(\xi')-\nabla_H\varphi(\xi')\cdot(\pm rv/\|v\|)\le K B_1 K_0 m^H_\varphi(\xi,r).
$$
\hfill$\square$\psn
\medskip

A result similar to Proposition \ref{prop tre in R}, involving now the $\Hn$-sections $\S_\varphi^\Hn (\xi,s)$, holds, but the proof is much more delicate:

\begin{proposition}\label{prop tre in R MAGARI}
Let $\varphi :\Hn \to \R$ be a function with round $H$-sections (with $K_0$ as in (\ref{round section})).
Then, there exists a constant $B_5>0$ such that, if $\xi'\in \S_\varphi^\Hn (\xi,s)$, then
\begin{equation}\label{B_5 in H}
 m^H_\varphi(\xi',r)\le B_5 m^H_\varphi(\xi,r),
\end{equation}
for $r$ such that  $s= m^H_\varphi(\xi,r)$. The constant $B_5$ depends only on $K_0$.
\end{proposition}

\noindent \textbf{Proof.} Since $\varphi$ has round $H$-sections, it belongs to
 $E(H,K)$, where $K$ depends only on $K_0$, and for every $r\ge 0$ and $\xi\in\Hn$ the inequality (\ref{B_4 in H zero}) holds.
 In addition, by Proposition \ref{proposizione proprieta m H},  the function $m_\varphi^H(\xi,\cdot):[0,+\infty)\to [0,+\infty)$ is invertible.

Take any $\xi_3\in  \S_\varphi^\Hn (\xi_0,s)$, i.e. $\xi_1\in  S_\varphi^H (\xi_0,s)$,  $\xi_2\in  S_\varphi^H (\xi_1,s)$ and  $\xi_3\in  S_\varphi^H (\xi_2,s)$, with $s$ such that $s=m^H_\varphi(\xi_0,r)$. By Proposition \ref{prop tre in R} and  $\xi_1\in  S_\varphi^H (\xi_0,s)$, we have
\begin{equation}\label{uno}
\left(m_\varphi^H(\xi_0,\cdot)\right)^{-1}(s)\le \left(m_\varphi^H(\xi_1,\cdot)\right)^{-1}(C_1s).
\end{equation}
Similarly,
 since $\xi_2\in  S_\varphi^H (\xi_1,s)$, we have
\begin{equation}\label{uno 1}
\left(m_\varphi^H(\xi_1,\cdot)\right)^{-1}(s)\le \left(m_\varphi^H(\xi_2,\cdot)\right)^{-1}(C_1s).
\end{equation}
Let us prove that there exists a constant $C$, which depends only on $C_1$ and $B_4,$ and hence on $K_0,$ such that
\begin{equation}\label{uno 2}
\left(m_\varphi^H(\xi,\cdot)\right)^{-1}(C_1s)
\le C
\left(m_\varphi^H(\xi,\cdot)\right)^{-1}(s),\qquad\forall \xi\in\Hn.
\end{equation}
Inequality (\ref{B_4 in H zero}) is equivalent to
$$
\left(m_\varphi^H(\xi,\cdot)\right)^{-1}(B_4\tilde s)\le 2\left(m_\varphi^H(\xi,\cdot)\right)^{-1}(\tilde s),\qquad \forall \tilde s\ge 0;
$$
by choosing $\beta\in\N$ such that $C_1\le B_4^\beta$, iterating the previous inequality and taking into account that $\tilde s\mapsto  \left(m_\varphi^H(\xi,\cdot)\right)^{-1}(\tilde s)$ is an increasing function we obtain
\begin{equation}\label{uno 3}
\left(m_\varphi^H(\xi,\cdot)\right)^{-1}(C_1s)
\le  2^\beta\left(m_\varphi^H(\xi,\cdot)\right)^{-1}\left(\frac{C_1s}{B_4^\beta}\right)
\le  2^\beta\left(m_\varphi^H(\xi,\cdot)\right)^{-1}(s).
\end{equation}
Hence, (\ref{uno 2}) holds with $C= 2^\beta$; now, by (\ref{uno}), (\ref{uno 3}) and (\ref{uno 1}), we obtain
\begin{eqnarray*}
\left(m_\varphi^H(\xi_0,\cdot)\right)^{-1}(s)
&\le& \left(m_\varphi^H(\xi_1,\cdot)\right)^{-1}(C_1s)\\
&\le&  2^\beta\left(m_\varphi^H(\xi_1,\cdot)\right)^{-1}(s)\\
&\le&  2^\beta\left(m_\varphi^H(\xi_2,\cdot)\right)^{-1}(C_1s).
\end{eqnarray*}
A similar argument proves that $\xi_3\in  S_\varphi^H (\xi_2,s)$ implies
$$
\left(m_\varphi^H(\xi_0,\cdot)\right)^{-1}(s)
\le 2^{2\beta}
\left(m_\varphi^H(\xi_3,\cdot)\right)^{-1}(C_1s).
$$
Now, recalling that $s= m^H_\varphi(\xi_0,r)$, the previous inequality gives
\begin{equation}\label{uno 2 bis}
m_\varphi^H(\xi_3,2^{-2\beta}r)\le C_1 m^H_\varphi(\xi_0,r).
\end{equation}
Finally, (\ref{B_2 in H zero}) and (\ref{uno 2 bis}) implies
$$
m_\varphi^H(\xi_3,r)\le B_2^{2\beta} C_1 m_\varphi^H(\xi_3,r)\le  B_2^{2\beta} C_1 m^H_\varphi(\xi_0,r),
$$
and the proof in finished.
\hfill$\square$\psn
\medskip

In order to introduce and prove the main result of the section, we need the following
\begin{lemma}\label{lemma}
  Let $\varphi :\Hn \to \R$ be
an $H$-convex function with round $H$-sections.
\begin{itemize}
\item[a.] If $\xi'\in  S_\varphi^H (\xi,m_\varphi^H(\xi,r))$ for some $r>0$, then
\begin{equation}\label{dis lemma}
 m_\varphi^H(\xi,r)\le m_\varphi^H(\xi',2^{1+\gamma}r),\end{equation}
 with  $\gamma$ as in (\ref{gamma});
 \item[b.] if $\xi'\in  S_\varphi^H (\xi, M_\varphi^H(\xi,r))$ for some $r>0$, then
\begin{equation}\label{dis lemma due}
 M_\varphi^H(\xi',r)\le m_\varphi^H(\xi,2^{2+3\gamma+\tilde \gamma}r),
\end{equation}
with $\tilde\gamma$ as in (\ref{tilde gamma}) which depends only on $K_0$.
\end{itemize}
\end{lemma}

\noindent \textbf{Proof.}
Let us consider $ \xi'=\xi\circ\exp v\in  S_\varphi^H (\xi,m_\varphi^H(\xi,r))$: by (\ref{inclusione con palla 2}) we have
\begin{equation}\label{dim new 1}\|v\|\le r.\end{equation}
The $H$-convexity of $\varphi$ and the $H$-monotonicity of $\nabla_H\varphi$ give
\begin{eqnarray}
  &&\varphi(\xi)-\varphi(\xi\circ\exp v)-\nabla_H\varphi(\xi\circ\exp v)\cdot v \ge 0\label{dim new 2}\\
 &&\left( \nabla_H\varphi(\xi\circ\exp v)-\nabla_H\varphi(\xi)\right)\cdot v \frac{r}{\|v\|}\ge 0;\label{dim new 3}
\end{eqnarray}
Again the $H$-convexity of $\varphi$ and (\ref{dim new 1})--(\ref{dim new 3}) give
\begin{eqnarray*}
 m_\varphi^H(\xi,r)&\le &
 \varphi\left(\xi\circ\exp\left(-\frac{v}{\|v\|}r\right)\right)-\varphi(\xi)-\nabla_H\varphi(\xi)\cdot\left(-\frac{v}{\|v\|}r\right)\\
 &\le &
 \varphi\left(\xi\circ\exp\left(-\frac{v}{\|v\|}r\right)\right)
 -\varphi(\xi\circ\exp v)- \nabla_H\varphi(\xi\circ \exp v)\cdot \left(-v \frac{r}{\|v\|}-v\right)
\\
&\le& M_\varphi^H(\xi',2r).
\end{eqnarray*}
Therefore (\ref{dis lemma}) follows from \eqref{remark_gamma}.

Let us prove b. Take any  $\xi'\in  S_\varphi^H (\xi, M_\varphi^H(\xi,r))$; since $\varphi\in E(H,K),$ with $K$ depending on $K_0$ only (see Theorem
\ref{round implica engulfing}), then
$\xi\in  S_\varphi^H (\xi',K M_\varphi^H(\xi,r))$. From \eqref{remark_gamma}
we have
%$ M_\varphi^H(\xi,r)\le  m_\varphi^H(\xi,2^\gamma r);$ hence
$ \xi'\in  S_\varphi^H (\xi, M_\varphi^H(\xi,r))\subset  S_\varphi^H (\xi, m_\varphi^H(\xi,2^\gamma r)),$ and (\ref{dis lemma}) implies that
\begin{equation}\label{dim lemma}
K M_\varphi^H(\xi,r)\le
K m_\varphi^H(\xi,2^\gamma r)\le
K m_\varphi^H(\xi',2^{2\gamma+1} r).
\end{equation}
Now, let $\tilde \gamma \in \N$ be such that
\begin{equation}\label{tilde gamma}
K\le B_4^{\tilde\gamma}.
\end{equation}
By iterating inequality (\ref{B_4 in H zero}) and  (\ref{remark_gamma}), inequality (\ref{dim lemma})  gives
\begin{equation}\label{dim lemma 1}
K M_\varphi^H(\xi,r)\le K m_\varphi^H(\xi',2^{2\gamma+1} r)\le \frac{K}{B_4^{\tilde\gamma}}m_\varphi^H(\xi',2^{1+2\gamma+\tilde\gamma}r)\le m_\varphi^H(\xi',2^{1+2\gamma+\tilde\gamma}r):
\end{equation}
Hence, $\xi\in  S_\varphi^H (\xi',K M_\varphi^H(\xi,r))\subset  S_\varphi^H (\xi',m_\varphi^H(\xi',2^{1+2\gamma+\tilde\gamma}r)).$ Finally, the inequalities
(\ref{remark_gamma}) and (\ref{dis lemma}) imply
$$
M_\varphi^H(\xi',r)\le  m_\varphi^H(\xi',2^{\gamma} r)\le m_\varphi^H(\xi',2^{1+2\gamma+\tilde\gamma}r)\le  m_\varphi^H(\xi,2^{2+3\gamma+\tilde\gamma}r).
$$
\hfill$\square$\psn
\medskip

We are now in the position to prove the first part of our main result in Theorem \ref{intro Teo}:

\noindent \textbf{Proof of Theorem \ref{intro Teo} i.}
Let $\varphi :\Hn \to \R$ be
an $H$-convex function with round $H$-sections. Let us prove that
 $\varphi$   satisfies the engulfing property $E(\Hn, K)$.
Fix $\xi\in\Hn$ and $s>0$.
Let us suppose that $\xi'\in \S_\varphi^\Hn(\xi,s)$: we have to prove that
 $\S_\varphi^\Hn(\xi,s)\subset  \S_\varphi^\Hn(\xi',Ks)$, where $ K$ is a  constant which depends only on $K_0$ in (\ref{round section}).

Let  $r$ be such that $s=m_\varphi^H(\xi,r)$. By definition,
\begin{equation}\label{dim new 0}
\S_\varphi^\Hn\left(\xi, s\right)=\bigcup_{
\begin{array}{l}
\xi_1\in S^H_\varphi(\xi,  m_\varphi^H(\xi,r))\\
 \xi_2\in S^H_\varphi(\xi_1,  m_\varphi^H(\xi,r))
 \end{array}
} \hskip -1.3truecm
S^H_\varphi(\xi_2,   m_\varphi^H(\xi,r)).
\end{equation}
For every $\xi_1\in   S_\varphi^H (\xi,m_\varphi^H(\xi,r)),$ by applying a. in Lemma \ref{lemma}, we get
\begin{equation}\label{dim new 4 bis}
 m_\varphi^H(\xi, r)\le m_\varphi^H(\xi_1,2^{1+\gamma}r)
\end{equation}
and, by (\ref{dim new 0}), we get
\begin{equation}\label{dim new 5}
\S_\varphi^\Hn\left(\xi, s\right)\subset\bigcup_{
\begin{array}{l}
\xi_1\in S^H_\varphi(\xi,  m_\varphi^H(\xi,r))\\
 \xi_2\in S^H_\varphi(\xi_1,   m_\varphi^H(\xi_1,2^{1+\gamma}r))
 \end{array}
} \hskip -1.7truecm
S^H_\varphi(\xi_2,   m_\varphi^H(\xi,r)).
\end{equation}
For every  $ \xi_2\in  S_\varphi^H (\xi_1,m_\varphi^H(\xi_1,2^{1+\gamma} r))$ with $\xi_1\in S^H_\varphi(\xi,  m_\varphi^H(\xi,r)),$ via  a. in Lemma \ref{lemma} we get
\begin{equation}\label{dim new 6}
 m_\varphi^H(\xi_1,2^{1+\gamma} r)\le m_\varphi^H(\xi_2,2^{2+2\gamma}r).
\end{equation}
Using now (\ref{dim new 4 bis}) and (\ref{dim new 6}), relation (\ref{dim new 5})
becomes
$$
\S_\varphi^\Hn\left(\xi, s\right)\subset\bigcup_{
\begin{array}{l}
\xi_1\in S^H_\varphi(\xi,  m_\varphi^H(\xi,r))\\
 \xi_2\in S^H_\varphi(\xi_1,   m_\varphi^H(\xi_1,2^{1+\gamma}r))
 \end{array}
} \hskip -1.7truecm
S^H_\varphi(\xi_2, m_\varphi^H(\xi_2,2^{2+2\gamma}r)  ).
$$
Since $\gamma>0$, we have
 the following inclusions:
\begin{eqnarray}
 \S_\varphi^\Hn\left(\xi, s\right)&\subset& \bigcup_{
\begin{array}{l}
\xi_1\in S^H_\varphi(\xi,  m_\varphi^H(\xi,2^{2+2\gamma}r))\\
 \xi_2\in S^H_\varphi(\xi_1,   m_\varphi^H(\xi_1,2^{2+2\gamma}r))
 \end{array}
} \hskip -1.7truecm
S^H_\varphi(\xi_2, m_\varphi^H(\xi_2,2^{2+2\gamma}r)  )\nonumber\\
&\subset&
\widetilde B(\xi,2^{2+2\gamma}r),\label{dim new 7}
\end{eqnarray}
where the last inclusion comes from  (\ref{gen incl 1}). Then, using the inclusions in (\ref{inclusione palle}), we get
\begin{equation}\label{dim new 8}
 \S_\varphi^\Hn\left(\xi, s\right)\subset\widetilde B(\xi,2^{2+2\gamma}r)\subset   B_g(\xi,3\,2^{2+2\gamma}r)\subset  B_g(\xi',3\,2^{3+2\gamma}r)
\subset   \widetilde B(\xi',3C2^{3+2\gamma}r),
\end{equation}
where  $C$  is the constant in the  Folland--Stein Lemma. Inclusions  (\ref{gen incl 2}) and (\ref{dim new 8}) give
$$
%\begin{equation}\label{dim new 9}
\S_\varphi^\Hn\left(\xi, s\right)\subset
  \bigcup_{
\begin{array}{l}
\xi_3\in S^H_\varphi(\xi',  M_\varphi^H(\xi',3C2^{3+2\gamma}r))\\
 \xi_4\in S^H_\varphi(\xi_3,  M_\varphi^H(\xi_3,3C2^{3+2\gamma}r))
 \end{array}
} \hskip -1.7truecm
 S^H_\varphi(\xi_4,   M_\varphi^H(\xi_4,3C2^{3+2\gamma}r)).
%\end{equation}
$$
Now, applying twice (\ref{dis lemma due}) in Lemma \ref{lemma}, we have, by the previous inclusion,
\begin{eqnarray}
\S_\varphi^\Hn\left(\xi, s\right)
&\subset&
  \bigcup_{
\begin{array}{l}
\xi_3\in S^H_\varphi(\xi',  M_\varphi^H(\xi',3C2^{3+2\gamma}r))\\
 \xi_4\in S^H_\varphi(\xi_3,  M_\varphi^H(\xi_3,3C2^{3+2\gamma}r))
 \end{array}
} \hskip -1.7truecm
 S^H_\varphi(\xi_4,   M_\varphi^H(\xi_3,3C2^{5+5\gamma+\tilde\gamma}r))\nonumber\\
 &\subset&
  \bigcup_{
\begin{array}{l}
\xi_3\in S^H_\varphi(\xi',  M_\varphi^H(\xi',3C2^{5+5\gamma+\tilde\gamma}r))\\
 \xi_4\in S^H_\varphi(\xi_3,  M_\varphi^H(\xi_3,3C2^{5+5\gamma+\tilde\gamma}r))
 \end{array}
} \hskip -1.7truecm
 S^H_\varphi(\xi_4,   M_\varphi^H(\xi_3,3C2^{5+5\gamma+\tilde\gamma}r))\nonumber\\
  &\subset&
  \bigcup_{
\begin{array}{l}
\xi_3\in S^H_\varphi(\xi',  M_\varphi^H(\xi',3C2^{5+5\gamma+\tilde\gamma}r))\\
 \xi_4\in S^H_\varphi(\xi_3,  M_\varphi^H(\xi',3C2^{7+8\gamma+2\tilde\gamma}r))
 \end{array}
} \hskip -1.7truecm
 S^H_\varphi(\xi_4,   M_\varphi^H(\xi',3C2^{7+8\gamma+2\tilde\gamma}r))\nonumber\\
  &\subset&
  \bigcup_{
\begin{array}{l}
\xi_3\in S^H_\varphi(\xi',  M_\varphi^H(\xi',3C2^{7+8\gamma+2\tilde\gamma}r))\\
 \xi_4\in S^H_\varphi(\xi_3,  M_\varphi^H(\xi',3C2^{7+8\gamma+2\tilde\gamma}r))
 \end{array}
} \hskip -1.7truecm
 S^H_\varphi(\xi_4,   M_\varphi^H(\xi',3C2^{7+8\gamma+2\tilde\gamma}r))\nonumber\\
 &=&\S_\varphi^\Hn\left(\xi',  M_\varphi^H(\xi',3C2^{7+8\gamma+2\tilde\gamma}r)\right)\label{dim new 10}.
\end{eqnarray}
Set $\tilde C=3C2^{7+8\gamma+2\tilde\gamma},$ and take any $\delta\in\N$ such that $\tilde C\le 2^\delta$; clearly, both $\tilde C$ and $\delta$ they depend only on $K_0$.
Hence, we have
 the following inclusions:
\begin{eqnarray}
 \S_\varphi^\Hn\left(\xi, s\right)
&\subset&  \S_\varphi^\Hn\left(\xi', M_\varphi^H\left(\xi',\tilde C r\right) \right) \nonumber\\
&\subset&  \S_\varphi^\Hn\left(\xi', B_1^\delta M_\varphi^H\left(\xi',\tilde C 2^{-\delta}r\right)\right )\qquad \qquad \qquad (\texttt{\rm by}\ \ref{B_1 in H zero}))\nonumber\\
&\subset&  \S_\varphi^\Hn\left(\xi', B_1^\delta  M_\varphi^H\left(\xi',r\right)\right)\nonumber\\
&\subset&  \S_\varphi^\Hn\left(\xi', K_1 B_1^\delta m_\varphi^H\left(\xi',r\right) \right)\qquad \qquad \qquad (\texttt{\rm by}\ (\ref{H mM})) \nonumber\\
&\subset&  \S_\varphi^\Hn\left(\xi', B_5 K_1 B_1^\delta m_\varphi^H\left(\xi,r\right)\right )\qquad \qquad \qquad (\texttt{\rm by}\ (\ref{B_5 in H})) \nonumber\\
&\subset&  \S_\varphi^\Hn\left(\xi' ,B_5 K_1 B_1^a s \right). \nonumber
\end{eqnarray}
\hfill$\square$\psn

\section{Balls and quasi-metrics via the $\Hn$-sections of $H$-convex functions}

It is known that there is a deep connection between the existence of a quasi-metric $d$ on a given set $X\subset\R^k,$ and the existence of a family of subsets  $\{S(x,s)\}_{\{x\in X,\, s>0\}}$ enjoying the following properties
\begin{itemize}
\item[$(P_1)$] $\bigcap_{s>0} S(x,s)=\{x\},$ for every $x\in X$;

\item[$(P_2)$] $\bigcup_{s>0} S(x,s)=\R^k,$ for every $x\in X$;

\item[$(P_3)$] for each $x\in X$, $s\mapsto S(x,s)$ is a non decreasing map;

\item[$(P_4)$] there exists a constant $H$ such that, for all $y\in S(x,s)$,
\begin{eqnarray}
S(x,s)&\subset& S(y,Hs),\label{prop 1}\\
S(y,s)&\subset& S(x,Hs).\label{prop 2}
\end{eqnarray}
\end{itemize}
As a matter of fact, the following result holds:
\begin{lemma}\label{lemma family}(see Lemma 1 in \cite{AiFoTo1998})
Let $X$ be a set and $S: X\times\R^+\to {\mathcal P}(X)$ be a set-valued map such that the family $\{S(x,s)\}$ has the properties $(P_1)$-$(P_4)$. Then, the function $d:X\times X\to [0,+\infty)$ defined by
$$
d(x,y)=\inf\left\{s:\ x\in S(y,s),\ y\in S(x,s)\right\}
$$
is a quasi-metric. On the other hand, given a quasi-metric $d$ defined on $X$, the family of the $d$-balls in $X$ satisfies the properties $(P_1)$-$(P_4)$.
\end{lemma}
In particular, in \cite{AiFoTo1998} the authors prove that the sections $S_u(x,r)$ of a convex function $u:\R^k\to\R$ satisfying the engulfing property, generate a quasi-metrics.

Let us now consider an $H$-convex function $\varphi:\Hn\to\R$ with round $H$-sections; by taking all $s>0$ and $\xi\in\Hn=\R^{2n+1}$ we obtain a family of sets $\{\S_\varphi^\H\left(\xi,s \right)\}_{\{\xi\in \Hn,\, s>0\}}$ (the $\Hn$-sections) for which conditions $(P_1)$-$(P_3)$ trivially hold; moreover, due to
Theorem \ref{intro Teo}, such family satisfies the engulfing property $E(\Hn,K)$, i.e. condition $(\ref{prop 1})$.

The next result shows that the family of $\Hn$-sections satisfies condition $(\ref{prop 2})$ too:
\begin{theorem}\label{round e topo}
Let $\varphi :\Hn \to \R$ be
an $H$-convex function with round $H$-sections. Then, there exists a constant $\tilde K,$ which depends only on $K_0,$ such that,
if $\xi'\in \S_\varphi^\Hn(\xi,s)$, then  $\S_\varphi^\Hn(\xi',s)\subset  \S_\varphi^\Hn(\xi,\tilde K s)$.
\end{theorem}

\noindent \textbf{Proof.} The proof follows the ideas in the proof of Theorem \ref{intro Teo}. Fix $\xi\in\Hn\ s>0$ and $\xi'\in \S_\varphi^\Hn(\xi,s)$; let  $r$ be such that $s=m_\varphi^H(\xi',r)$.
Theorem  \ref{intro Teo}
guarantees that
$\varphi$   satisfies the engulfing property $E(\Hn, K)$, where $K$ depends only on $K_0$.
Hence, $\xi\in \S_\varphi^\Hn(\xi',Ks)$.
Proposition \ref{prop tre in R MAGARI} implies that
$$
 m^H_\varphi(\xi,\widehat r)\le B_5 m^H_\varphi(\xi',\widehat r),
$$
for $\widehat r$ such that  $Ks= m^H_\varphi(\xi',\widehat r)$ (the constant $B_5$ depends only on $K_0$).
Since, by Proposition \ref{proposizione proprieta m H},  the function $r\mapsto m^H_\varphi(\xi, r)$ is an increasing function, we obtain
\begin{equation}\label{B_5 in H new}
 m^H_\varphi(\xi, r)\le B_5 m^H_\varphi(\xi',\widehat r)=B_5 Ks.
\end{equation}
By definition,
\begin{equation}\label{dim new 0 due}
\S_\varphi^\Hn\left(\xi', s\right)=\bigcup_{
\begin{array}{l}
\xi_1\in S^H_\varphi(\xi',  m_\varphi^H(\xi',r))\\
 \xi_2\in S^H_\varphi(\xi_1,  m_\varphi^H(\xi',r))
 \end{array}
} \hskip -1.4truecm
S^H_\varphi(\xi_2,   m_\varphi^H(\xi',r)).
\end{equation}
Using exactly the same arguments as in the proof of Theorem  \ref{intro Teo}, that allow us to pass from (\ref{dim new 0}) to
(\ref{dim new 7}) (essentially, by exchanging the role of $\xi$ and $\xi'$), we obtain
\begin{equation}\label{dim new 7 due}
 \S_\varphi^\Hn\left(\xi', s\right)\subset
\widetilde B(\xi',2^{2+2\gamma}r)\subset
B_g(\xi',3\,2^{2+2\gamma}r).\end{equation}
Now, taking into account the definition of $\tilde\gamma$ in (\ref{tilde gamma}) and
 iterating inequality (\ref{B_4 in H zero}), we get
\begin{equation}\label{dim lemma 1 due}
Ks=K m_\varphi^H(\xi', r)\le \frac{K}{B_4^{\tilde\gamma}}m_\varphi^H(\xi',2^{\tilde\gamma}r)\le m_\varphi^H(\xi',2^{\tilde\gamma}r).
\end{equation}
Since $\xi\in \S_\varphi^\Hn(\xi',Ks)$, we obtain
\begin{equation}\label{dim new 0 tre}
\S_\varphi^\Hn\left(\xi', Ks\right)\subset\bigcup_{
\begin{array}{l}
\xi_1\in S^H_\varphi(\xi',  m_\varphi^H(\xi',2^{\tilde\gamma}r))\\
 \xi_2\in S^H_\varphi(\xi_1,  m_\varphi^H(\xi',2^{\tilde\gamma}r))
 \end{array}
} \hskip -1.7truecm
S^H_\varphi(\xi_2,   m_\varphi^H(\xi',2^{\tilde\gamma}r)).
\end{equation}
Using exactly the same arguments that allow us to pass from (\ref{dim new 0 due}) to
(\ref{dim new 7 due}) (essentially, by exchanging the role of $r$ with $2^{\tilde\gamma}r$), we obtain
\begin{equation}\label{dim new 7 tre}
\xi\in
\widetilde B(\xi',2^{2+2\gamma+\tilde \gamma}r)\subset
B_g(\xi',3\,2^{2+2\gamma+\tilde\gamma}r).\end{equation}
Now, taking into account that $\tilde\gamma>0$, relations (\ref{dim new 7 due}) and (\ref{dim new 7 tre}) give
\begin{equation}\label{dim new 8 due}
 \S_\varphi^\Hn\left(\xi', s\right)\subset
B_g(\xi',3\,2^{2+2\gamma}r)\subset
B_g(\xi,3\,2^{3+2\gamma+\tilde\gamma}r)
\subset   \widetilde B(\xi,3C2^{3+2\gamma+\tilde\gamma}r),
\end{equation}
where  $C$ in the previous inclusions is the constant in the Folland--Stein Lemma.
Using the same arguments that allow us to pass from (\ref{dim new 8}) to (\ref{dim new 10}) (essentially, by replacing $3C2^{3+2\gamma}r$ with $3C2^{3+2\gamma+\tilde\gamma}r$), we obtain
$$
\S_\varphi^\Hn\left(\xi', s\right)\subset
\S_\varphi^\Hn\left(\xi,  M_\varphi^H(\xi,3C2^{7+8\gamma+3\tilde\gamma}r)\right).$$
Set $\widehat C=3C2^{7+8\gamma+3\tilde\gamma}$ and take $\widehat\delta\in\N$ such that $\widehat C\le 2^{\widehat \delta}$; clearly, $\widehat C$ and $\widehat \delta$  depend only on $K_0$.
We have
 the following inclusions:
\begin{eqnarray}
 \S_\varphi^\Hn\left(\xi', s\right)
&\subset&  \S_\varphi^\Hn\left(\xi, M_\varphi^H\left(\xi,\widehat C r\right) \right) \nonumber\\
&\subset&  \S_\varphi^\Hn\left(\xi, B_1^{\widehat \delta} M_\varphi^H\left(\xi,\widehat C 2^{-\widehat\delta}r\right)\right )\qquad \qquad \qquad (\texttt{\rm by}\ \ref{B_1 in H zero}))\nonumber\\
&\subset&  \S_\varphi^\Hn\left(\xi, B_1^{\widehat \delta}  M_\varphi^H\left(\xi,r\right)\right)\nonumber\\
&\subset&  \S_\varphi^\Hn\left(\xi, K_1 B_1^{\widehat \delta} m_\varphi^H\left(\xi,r\right) \right)\qquad \qquad \qquad (\texttt{\rm by}\ (\ref{H mM})) \nonumber\\
&\subset&  \S_\varphi^\Hn\left(\xi, B_5 K K_1 B_1^{\widehat \delta} s\right )\qquad \qquad \qquad (\texttt{\rm by}\ (\ref{B_5 in H new})) \nonumber
\end{eqnarray}
which concludes the proof.
\hfill$\square$\psn
\medskip

We are now in the position to prove the second part of our main result in Theorem \ref{intro Teo}:

\noindent \textbf{Proof of Theorem \ref{intro Teo} ii.} Let $\varphi :\Hn \to \R$ be an $H$-convex function with round $H$-sections.
The previous arguments, together with Lemma \ref{lemma family} and Lemma 2 in \cite{AiFoTo1998}, give that
$$
d_\varphi(\xi,\xi')=\inf\left\{s>0:\ \xi\in \S_\varphi^\Hn(\xi',s),\  \xi'\in \S_\varphi^\Hn(\xi,s)\right\}
$$
is a quasi-metric in $\Hn$. Moreover, if $B_\varphi(\xi,r)$ denotes the $d_\varphi$-ball of center $\xi\in\Hn$ and radius $r>0$,
we have that there exists $H$ which depends only on $K_0$ in (\ref{round section}) such that
\begin{equation}\label{dim finale balle}
  \S_\varphi^\Hn\left(\xi, \frac{r}{2H}\right)\subset B_\varphi(\xi,r)\subset \S_\varphi^\Hn\left(\xi, r\right).
\end{equation}
\hfill$\square$\psn
\medskip

The definition of $\H^n$-sections via subsequent constructions of $H$-sections makes hard its description in terms of functional inequalities.
However, in  the very simple case of the function $\varphi:\H\to\R$ defined by $\varphi(x,y,t)=x^2+y^2,$ we are able to fully describe the set $\S_\varphi^\H(e, r)$ by providing explicitly the equation of its boundary. While in the Euclidean
 case the function $u(x)=\|x\|,$ with $x\in\R^n$, gives rise to the sections $S_u(x_0,s)=B^{\R^n}(x_0,\sqrt{s}),$ i.e., the usual balls in $\R^n,$ in the case of the first Heisenberg group $\H$ and with the mentioned function $\varphi$ we obtain $S_\varphi^\H(\xi_0,s)=\widetilde B(\xi_0,\sqrt{s}),$ and the family of $\H$-sections of $\varphi$ consists of the $\widetilde B$-balls in (\ref{balls tilde}).

\begin{example}\label{esempio}{\rm
Let us consider $\varphi:\H\to \R$ defined by $\varphi(x,y,t)=x^2+y^2$. This function is $\R^3$-convex, and hence $H$-convex. Since $\partial_H\varphi(x,y,t)=\{2(x,y)\}$, the horizontal section $S^H_\varphi(\xi_0,s)$ is given by
$$
S^H_\varphi(\xi_0,s)=\{\xi=(x,y,t)\in H_{\xi_0}:\ (x-x_0)^2+(y-y_0)^2<s\},
$$
for $\xi_0=(x_0,y_0,t_0)$ and $s>0$. Hence, for this  particular $\varphi,$ we have that
\begin{equation}\label{example 1}
S^H_\varphi(\xi_0,s)=B_g(\xi_0,\sqrt s)\cap H_{\xi_0},
\end{equation}
and, therefore,
$$
\overline{\S^\H_\varphi(\xi_0,s)}=
\widetilde B(\xi_0,\sqrt s).
$$
Since, from the definition of $\H$-section, $\S^\H_\varphi(\xi_0,s)=\xi_0\circ \S^\H_\varphi(e,s)$, we will focus on the particular case $\xi_0=e.$
We  claim that, for every $r>0$,
\begin{equation}\label{example 1 bis}
\overline{\S^\H_\varphi(e,r)}=
\widetilde B(e,\sqrt r)=\left\{\xi=(x,y,t):\ |t|\le \sqrt{3r+2\|(x,y)\|\sqrt r-\|(x,y)\|^2}\left(\sqrt r+\|(x,y)\|\right) \right\}.
\end{equation}
Let us try to give the idea of its construction. Fix $r>0$. First of all, note that
\begin{itemize}
\item $\widetilde B(e,\sqrt r)$ is radial with respect to the $t$-axis;
\item $\widetilde B(e,\sqrt r)$ is symmetric with respect to the $xy$-plane.
\end{itemize}
In particular, it is sufficient to identify the points of the set $\partial  \widetilde B(e,\sqrt r)$ in $\H\cap\{t\ge 0\}.$
To this purpose, let us consider the points
\begin{equation}\label{example 2}
\eta^\theta=\left(\sqrt r,0,0\right)\circ\left(\sqrt r\cos \theta,\sqrt r\sin\theta,0\right)\circ \left(\sqrt r\cos(2 \theta),\sqrt r\sin(2\theta),0\right),\qquad \texttt{\rm for}\ \theta\in[-2\pi/3,0].\end{equation}
Trivially, $\eta^0=(3\sqrt r,0,0)\in\partial \widetilde B(e,\sqrt r)$.
Let us motivate our choice in (\ref{example 2}). Let $v_i\in V_1\cong\R^2$, for $i=1,2,3$, and consider the point
\begin{equation}\label{example 2 bis}
\eta=(x,y,t)=\exp(v_1)\circ\exp(v_2)\circ\exp(v_3);
\end{equation}
we have  $(x,y)=v_1+v_2+v_3,$ and   $|t|/4$ is equal to the area of the polygon $P={\mathrm{co}}\{ (0,0),\ v_1,\ v_1+v_2,\ v_1+v_2+v_3\}\subset\R^2,$ where \lq\lq co\rq\rq\ denotes the convex hull (for details on this application of Stokes' Theorem, see, for example, Section 2.3 in \cite{CaDaPaTy2007}).
 In order to construct $\partial  \widetilde B(e,\sqrt r)\cap\{(x,y,t)\in\H:\ t\ge 0\}$ we restrict our attention to the points $\eta$ in (\ref{example 2 bis}) with the following features:
 \begin{itemize}
 \item   $\|v_i\|=\sqrt r;$
 \item the angles $\widehat{v_1,v_2}$ and   $\widehat{v_2,v_3}$ are equal to $\theta$
 \end{itemize}
 (this choice will be explained later on).
Due to the symmetries of $\widetilde B(e,\sqrt r)$, we set
\begin{equation}\label{example 2 tris}
v_1=\sqrt{r}(1,0),\quad  v_2=\sqrt{r}(\cos \theta,\sin\theta),\quad
v_3=\sqrt{r}(\cos (2 \theta),\sin(2\theta)).
\end{equation}
With this choice, from (\ref{example 2}) one simply gets that
\begin{equation}\label{example 3}
\eta^\theta=(x(\theta),y(\theta),t(\theta))=\left(\sqrt r \left(1+\cos \theta+\cos(2 \theta)\right),
\sqrt r \left(\sin \theta+\sin(2 \theta)\right),
-4r\sin\theta(1+\cos\theta)\right).
\end{equation}
Clearly, $t(\theta)\ge 0$ for $\theta\in[-\pi,0]$.
In the case $\theta=-2\pi/3$, $P$ turns out to be an equilateral triangle, and
$\eta^{-2\pi/3}=(0,0,\sqrt 3 r)$; in the case $\theta\in [-\pi,-2\pi/3)$, we have that $\eta^{\theta}$ is an interior point of $\widetilde B(e,\sqrt r)$. Therefore, we restrict our attention
to the points $\eta^\theta$ as in (\ref{example 2}).
Simple computations give that, for $\theta\in[-2\pi/3,0],$
\begin{eqnarray*}
&&d(\theta):=\|(x(\theta),y(\theta))\|=\sqrt r (1+2\cos\theta)\\
&&t(\theta)=4r\sqrt{1-\cos \theta}(1+\cos\theta).
\end{eqnarray*}
Note that, if $\theta=-\pi/3,$ the function $t(\theta)$ reaches its maximum $ 3\sqrt 3 r$ and, in this case, $d(-\pi/3)=2\sqrt r$.
Consider the change of variable $z=\sqrt r (1+2\cos\theta);$  due to the symmetry of  $\widetilde B(e,\sqrt r),$ we obtain that
$$
\left(z,0,\sqrt{3  r+2z\sqrt r-z^2}\left(\sqrt r+z\right)\right)\in \partial  \widetilde B(e,\sqrt r), \qquad \texttt{\rm for}\ z\in[0,3\sqrt r],
$$
and thus we get the expression in (\ref{example 1 bis}).

\begin{center}
    \includegraphics[width=4 cm]{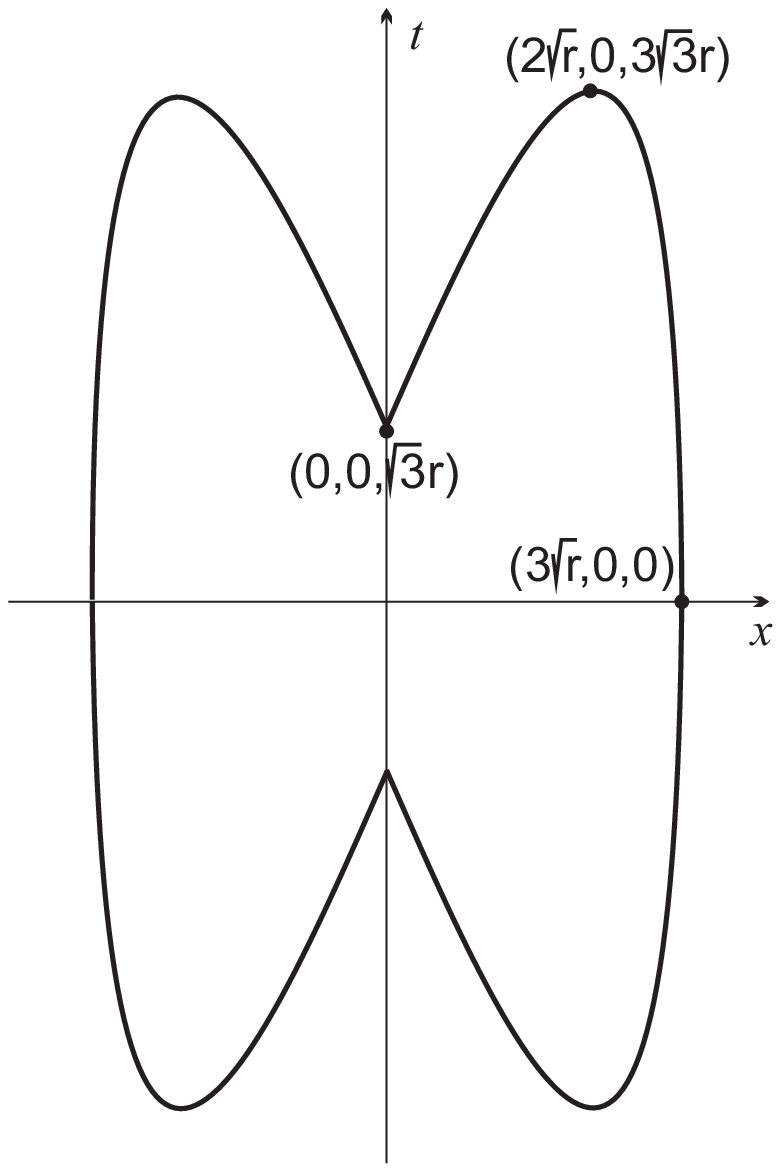}\\
    \emph{The profile in the plane $(x,0,t)$ of the $\H$-section $\S^\H_\varphi(e,r)$ of the function  $\varphi:\H\to \R$  defined by $\varphi(x,y,t)=x^2+y^2$, for $r=1$.}
         \end{center}
Finally, let us explain briefly the restrictions imposed in (\ref{example 2 tris}) to obtain (\ref{example 1 bis}).

First, it is easy to see that, if in (\ref{example 2 bis}) we set $\|v_i\|=\sqrt r'$, with $0<r'<r$, we obtain that $\eta$ in (\ref{example 2 bis}) is in $\partial\widetilde B(e,\sqrt{r'})\subset \widetilde B(e,\sqrt{r})$; a similar argument holds for $\eta$ in (\ref{example 2 bis}), with the choice $\|v_i\|<\sqrt r$.

Secondly, let us motivate the restriction $\widehat{v_1,v_2}=\widehat{v_2,v_3}=\theta$ in (\ref{example 3}). Fix $\theta\in  (-2\pi/3,0)$, consider $v_i$ as in (\ref{example 2 tris}) and the mentioned polygon $P$; using  (\ref{example 3}), the area of $P$ is exactly $-\sin\theta(1+\cos \theta).$ If one looks for the triplet of vectors  $v_i$, with $\|v_i\|=\sqrt r$ for $i=1,2,3,$  such that $v_1+v_2+v_3=(x(\theta),y(\theta))$ and  such that the area of the associated polygon $P$ is the biggest one, then one obtains exactly the vectors $v_i$ in (\ref{example 2 tris}). This proves that $\eta^\theta$ belongs to the boundary of our $\H$-section. We leave the details and their tedious calculations to the interested reader.}
\end{example}

\section{Final remarks and open questions}

\noindent \emph{Question 1.} The assumption of the round $H$-sections property for an $H$-convex function $\varphi$ is a sufficient condition in order to guarantee that $\varphi$ satisfies the engulfing property $E(\H^n,K).$ It would be nice to weaken this assumption and prove that a function with the  engulfing property $E(H,K)$ satisfies the engulfing property $E(\Hn,K)$.

\noindent \emph{Question 2.} In \cite{CapMal2006} the authors study the engulfing property for convex functions in a generic Carnot group $\textbf{G}$; as a matter of fact, in this more general framework, the related definition of $\textbf{G}$-sections (as in Definition \ref{definizione section Heis}) would be affected by the different geometry of the group $\textbf{G}$, by the number of the steps and, especially, by the number of consecutive horizontal segments needed to joint any pair of points. Moreover, in a Carnot group with step greater than 2, a so-called horizontal line, i.e., a set $\{\xi\circ \exp sv\}_{s\in \R},$ is not a line in the Euclidean sense, as well as a horizontal plane is not a hyperplane in the Euclidean sense. This leads us to think that the $\textbf{G}$-sections may have a very peculiar shape.

\noindent \emph{Question 3.}
In \cite{KoMa2005} the authors prove, among other things, that the notion of round sections in Definition \ref{def round sections}, controlled slope in (\ref{def controlled slope}), quasi uniform convexity, and quasiconformity are strictly related properties. To be precise, the next result holds (see Theorem 3.1 in \cite{KoMa2005}):
\begin{theorem}\label{teo kovalev maldonado} Let $n\ge 2$, and let $u:\R^n\to\R$ be a convex function. The following are equivalent:
\begin{itemize}
\item[i.] $u$ is quasiuniformly convex function, i.e. $u$ is not affine, $u\in W_{loc}^{2,n}$ and there exists a constant $K\ge 1$ such that
\begin{equation}\label{quasiuniformly convex}
\|\nabla^2 u (x)\|^n\le K {\rm det} \nabla^2 u(x),\qquad \mathrm{a.e.}\, x\in\R^n;
\end{equation}

\item[ii.] $u$ is differentiable and $\nabla u:\R^n\to\R^n$ is quasiconformal, recalling that an injective map $F:\R^n\to\R^n$ is quasiconformal if $F\in W_{loc}^{1,n}$ and there exists  a constant $K\ge 1$ such that
\begin{equation}\label{quasiconformal map}
\|\nabla F (x)\|^n\le K {\rm det} \nabla F (x),\qquad \mathrm{a.e.}\, x\in\R^n;
\end{equation}

\item[iii.] $u$ is differentiable, but not affine, and has controlled slope;
\item[iv.] $u$ has round sections.
\end{itemize}
\end{theorem}
On the other hand, it is well known that the notion of quasiconformal maps on $\Hn$ has been introduced and intensively studied (see for example \cite{CaDaPaTy2007}). In this paper we introduce the notion of $H$-controlled slope and round $H$-sections for an $H$-convex function but, at least  to our knowledge, a horizontal notion of quasiuniform convexity for $H$-convex function does not exist of the literature. Our future aim will be to investigate a horizontal version of Theorem \ref{teo kovalev maldonado}.

\end{document}